\documentclass[12pt,leqno]{article}
\usepackage{amsmath,amssymb,bbm,array}
\usepackage{amsthm}
\usepackage[all,cmtip]{xy}
\def\bm#1{\mathbbm{#1}}

\def\vds{\hbox{\hbox to 2\arraycolsep{\hss\vbox{\vbox to 1.3ex{ 
\vss\hbox{.}\vspace{.33ex}\hbox{.}\vspace{.33ex}\hbox{.}\vspace{.33ex}\hbox{.}\vspace{.33ex}\hbox{.}\vss}}\hss}}}

\makeatletter
\renewcommand{\section}{\@startsection{section}{1}{0mm}{5mm}{3mm}{\raggedright\bf\large}}

\def\@citex[#1]#2{\if@filesw\immediate\write\@auxout{\string\citation{#2}}\fi
  \def\@citea{}\@cite{\@for\@citeb:=#2\do
    {\@citea\def\@citea{\@citesep}\@ifundefined
       {b@\@citeb}{{\bf ?}\@warning
       {Citation `\@citeb' on page \thepage \space undefined}}%
{\csname b@\@citeb\endcsname}}}{#1}}
\def\@citesep{; }
\makeatother

\newtheoremstyle{Kang}{}{}{\itshape}{}{\bf}{}{.5em}{}
\theoremstyle{Kang}
\newtheorem{theorem}{Theorem}[section]

\newtheorem{proposition}[theorem]{Proposition}

\newtheorem{lemma}[theorem]{Lemma}

\newtheoremstyle{Kremark}{}{}{}{}{\bf}{}{.5em}{}
\theoremstyle{Kremark}

\newtheorem*{remark}{Remark.}
\newtheorem{other}{}

\allowdisplaybreaks[1]  
\numberwithin{equation}{section}

\begin{document}

\title{Rationality problem for transitive subgroups of $S_8$}
\author{
Baoshan Wang \; and \; Jian Zhou \\ \normalsize School of
Mathematics
and System Sciences, Beihang University, Beijing, China \\
\normalsize E-mail: bwang@buaa.edu.cn \\ \normalsize School of
Mathematical Sciences, Peking University, Beijing, China
\\ \normalsize E-mail: zhjn@math.pku.edu.cn }
\date{}

\maketitle


\setlength{\parskip}{3pt}

\begin{abstract}
{\bf Abstract.} For any field $K$ and any transitive subgroup $G$ of
$S_8$, let $G$ acts naturally on $K(x_1,\dots,x_8)$ by permutations
of the variables, we prove that under some minor conditions
$K(x_1,\dots,x_8)^G$ is always $K$-rational except $G$ is $A_8$ or
$G$ is isomorphic to ${\rm PGL}(2,7)$. We pay special attentions on
the characteristic 2 cases.
\end{abstract}

{\small \hspace*{4mm}2010 Mathematics Subject Classification.
Primary 12F20, 13A50, 14E08.

\hspace*{4mm}{\it Keywords}: rationality, symmetric group, Noether's problem.\\
}

\section{Introduction}

Let $K$ be a field, $x_1,\dots,x_n$ be variables, and $G$ be any
transitive subgroup of $S_n$. Then $G$ acts on the rational function
field $K(x_1,\dots,x_n)$ by permutations of the variables. In
Noether's approach to the inverse Galois problem (cf. \cite{No1},
\cite{No2}), it is important to know if the fixed field
$K(x_1,\dots,x_n)^G$ is $K$-rational (= a purely transcendental
extension of $K$) or not. When $n \le 7$, fairly complete results
have been achieved by \cite{KW}, \cite{KWZ}, \cite{Zh}. Recall first
that:

\begin{theorem}
Let $K$ be any field, $G$ be a subgroup of $S_n$ not equal to $A_n$.
Let $G$ acts on the rational function field $K(x_1,\dots,x_n)$ via
$K$-automorphisms defined by $\sigma(x_i) = x_{\sigma(i)}$ for any
$\sigma \in G$, any $1 \le i \le n$.

{\rm (1) (\cite{KW}, Th. 1.3)} If $n \le 5$, then
$K(x_1,\dots,x_n)^G$ is $K$-rational. {\rm (In this case, $G = A_n$
is also OK, while $A_3, A_4$ are easy, $A_5$ is due to {\rm
\cite{Mae}}.)}

{\rm (2) (\cite{KWZ}, Th. 1.2)} If $n = 6$, then
$K(x_1,\dots,x_6)^G$ is $K$-rational except $G$ is isomorphic to
${\rm PSL}(2,5)$ or ${\rm PGL}(2,5)$. When $G$ is isomorphic to
${\rm PSL}(2,5)$ or ${\rm PGL}(2,5)$, then $K(x_1,\dots,x_6)^G$ is
stably $K$-rational, and $\mathbb{C}(x_1,\dots,x_6)^G$ is
$\mathbb{C}$-rational.

{\rm (3) (\cite{KW}, Th. 1.4)} If $n = 7$, and $G$ is a transitive
subgroup of $S_7$, then $K(x_1,\dots,$ $x_7)^G$ is $K$-rational
except $G$ is isomorphic to ${\rm PSL}(2,7)$. When $G$ is isomorphic
to ${\rm PSL}(2,7)$, and $K$ contains $\mathbb{Q}(\sqrt{-7})$, then
$K(x_1,\dots,x_7)^G$ is $K$-rational.

{\rm (4) (\cite{KW}, Th. 1.5)} If $n = 11$, and $G$ is a transitive
solvable subgroup of $S_{11}$, then $K(x_1,\dots,x_{11})^G$ is
$K$-rational.
\end{theorem}

In this paper, we will continue these investigations by considering
the case $n=8$, our main result is

\begin{theorem} Let $K$ be any field, $G$ be a transitive subgroup
of $S_8$ not equal to $A_8$, and consider the $K$-linear action of
$G$ on $K(x_1,\dots,x_8)$ by $g(x_i) = x_{g(i)}$ for all $g \in G$
and $1 \le i\le 8$.

{\rm (1)} For $G$ solvable, the fixed field $K(x_1,\dots,x_8)^G$ is
$K$-rational except the case when $G$ is conjugate to $\langle
(1,2,3,4,5,6,7,8)\rangle \simeq C_8$. For this last group,
$K(x_1,\dots,x_8)^G$ is $K$-rational if and only if ${\rm char}(K) =
2$ or ${\rm Gal}(K(\zeta_8)/K)$ is cyclic, where $\zeta_8$ is a
primitive 8th root of unity.

{\rm (2)} For $G$ unsolvable, assume in addition that $K$ contains
$\mathbb{Q}(\sqrt{-7})$. Then $K(x_1,\dots,$ $x_8)^G$ is
$K$-rational except when $G$ is isomorphic to ${\rm PGL}(2,7)$.
\end{theorem}

The rationality for $A_n$ when $n \ge 6$ is a famous unsolved
problem, while the rationality for ${\rm PGL}(2,7)$ seems to be out
of reach by the techniques of the present paper.

The paper is organized in the following way. In \S2, we collect many
results that frequently used in the rationality proofs. In \S3, we
present a list of transitive subgroups of $S_8$ based on \cite{CHM},
and we divide these groups (except one) into 4 types. Then in the
following \S4, \S5, \S6, \S7, we give the rationality proofs for the
4 types separately.

Many computations in this paper are done with the help of Wolfram
Mathematica.


\section{Preliminaries}

Let $L$ be an finitely generated extension field of $K$, with
transcendence degree $n$. To prove that $L$ is $K$-rational, it is
enough to show that $L$ can be generated by $n$ elements over $K$.
For $n \ge 3$, our strategy of rationality proof is to find out
these $n$ elements. By using the following basic theorem, we can
reduce our problem to a lower dimensional one.

\begin{theorem} Let $G$ be a finite group acting on
$L(x_1,\dots,x_n)$, the rational function field over a field $L$ on
$n$ variables $x_1,\dots,x_n$.

{\rm (1) (\cite{HaK3}, Th. 1)} Assume that $L$ is stable under $G$,
and the restricted action of $G$ on $L$ is faithful. Assume in
addition that $G$ acts on $x_i$'s affine-linearly, i.e. for any
$\sigma \in G$, we have
$$
\begin{pmatrix}
\sigma(x_1) \\ \vdots \\ \sigma(x_n)
\end{pmatrix} = A(\sigma) \begin{pmatrix}
x_1 \\ \vdots \\ x_n
\end{pmatrix} + B(\sigma), $$ where $A(\sigma) \in {\rm GL}_n(L)$ and
$B(\sigma) \in L^n$. Then $L(x_1,\dots,x_n) = L(z_1,\dots,z_n)$ for
some $z_1,\dots,z_n$ that are invariant under $G$. In particular,
$L(x_1,\dots,x_n)^G = L^G(z_1,\dots,z_n)$ is rational over $L^G$.

{\rm (2) (\cite{Mi}, see also \cite{AHK}, Th. 3.1)} In the above
situation, if $n = 1$, the same result holds without the assumption
that $G$ acts on $L$ faithfully.
\end{theorem}

A concrete form of this theorem for permutation groups is given by:

\begin{proposition} Let $S_n$ acts on $x_1,\dots,x_n$ and
$x_1^{(j)},\dots,x_n^{(j)}\; (1 \le j \le m)$ by permuting the lower
indexes, namely $\sigma(x_i) = x_{\sigma(i)}, \sigma(x_i^{(j)}) =
x_{\sigma(i)}^{(j)}$ for $\sigma \in S_n$. Let $K$ be any field, and
$G$ be any subgroup of $S_n$, which acts on
$x_1^{(j)},\dots,x_n^{(j)}\; (1 \le j \le m)$ by restriction, then
we have
\begin{multline*}
K(x_1,\dots,x_n,x_1^{(1)},\dots,x_n^{(1)},\dots,x_1^{(m)},\dots,x_n^{(m)})^G
\\ =
K(x_1,\dots,x_n)^G(t_1^{(1)},\dots,t_n^{(1)},\dots,t_1^{(m)},\dots,t_n^{(m)}),
\end{multline*} where $$
\begin{array}{ccccccc} t_1^{(j)} &=& x_1^{(j)} &+& \cdots &+& x_n^{(j)}, \\
t_2^{(j)} &=& x_1 x_1^{(j)} &+& \cdots &+& x_n x_n^{(j)}, \\ &\vdots& \\
t_n^{(j)} &=& x_1^{n-1} x_1^{(j)} &+& \cdots &+& x_n^{n-1} x_n^{(j)}
\end{array}
$$ for $1 \le j \le m$.
\end{proposition}

The rationality problem for a wreath product of permutation groups
can be reduced to corresponding problems for each factor.

\begin{theorem} {\rm ({\rm \cite{KWZ}, Th. 2.5})}
Let $K$ be any field, $G \subseteq S_m$ and $H \subseteq S_n$. Let
$G$ and $H$ act on the rational function fields $K(x_1, \dots ,
x_m)$ and $K(y_1, \dots , y_n)$ respectively via $K$-automorphisms
defined by $g(x_i) = x_{g(i)},\, h(y_j) = y_{h(j)}$ for any $g \in
G, h \in H, 1 \le i \le m, 1 \le j \le n$. Then $T := H \wr G$ may
be regarded as a subgroup of $S_{mn}$ acting on the rational
function field $K((x_{ij})_{1 \le i \le m, 1 \le j \le n})$. Assume
that both $K(x_1, \dots , x_m)^G$ and $K(y_1, \dots , y_n)^H$ are
$K$-rational. Then $K((x_{ij})_{1 \le i \le m, 1 \le j \le n})^T$ is
also $K$-rational.
\end{theorem}

In many cases, rationality problems for permutation actions can be
reduced to corresponding problems for monomial actions. Recall that
an $K$-action of a finite group $G$ on the rational function field
$K(x_1,\dots,x_n)$ is called monomial, if for any $g \in G$ we have
$$ g(x_j) = c_j(g) \prod_{i=1}^n x_i^{a_{ij}(g)}
$$ for $(a_{ij}(g)) \in {\rm GL}(n,\mathbb{Z})$ and $c_j(g) \in K$
(necessarily nonzero). When all $c_j(g) = 1$, the action is called
purely monomial.

\begin{theorem} Let $G$ be a finite group, $K(x_1,\dots,x_n)$ be the
rational function field endowed with a monomial $K$-action of $G$.

{\rm (1) (\cite{Ha1}, \cite{Ha2})} If $n = 2$, the fixed field
$K(x_1,x_2)^G$ is $K$-rational.

{\rm (2) (\cite{HaK1}, \cite{HaK2}, \cite{HoRi})} If $n = 3$, and
the action of $G$ is purely monomial, then $K(x_1,x_2,x_3)^G$ is
$K$-rational.
\end{theorem}

A monomial action has an associated integral representation $G \to
{\rm GL}(n,\mathbb{Z})$, namely $g \mapsto (a_{ij}(g))$. We say a
monomial action is reduced, if the associated integral
representation is faithful. By Lem. 2.8 of \cite{KP}, rationality
problems for monomial actions can always be reduced to those of the
reduced actions.

\begin{theorem} {\rm (\cite{HoKiYa}, Th. 1.6)} Let $K$ be a field of characteristic
$\ne 2$, $G$ be a finite group, and $K(x_1,x_2,x_3)$ be the rational
function field on three variables $x_1,x_2,x_3$ endowed with a
reduced monomial $K$-action of $G$. Then $K(x_1,x_2,x_3)^G$ is
$K$-rational except when the image of $G$ in ${\rm
GL}(3,\mathbb{Z})$ is conjugate to the following 9 groups $$
G_{1,2,1},\; G_{2,3,1},\; G_{3,1,1},\; G_{3,3,1},\; G_{4,2,1},\;
G_{4,2,2},\; G_{4,3,1},\; G_{4,4,1},\; G_{7,1,1}. $$ {\rm (See
\cite{HoKiYa}, section 2 for the definitions of these groups.)}
\end{theorem}

\cite{HoKiYa} also contains some partial results on $G_{7,1,1}$. For
the other 8 groups, \cite{Ya} gives necessary and sufficient
conditions on the field $K$ for the fixed field to be $K$-rational.

\bigskip

Sometimes homogeneity method is used in a rationality proof. Here is
a brief discussion.

Let $L$ be an extension field of $K$, and assume that $L$ has a
sub-$K$-algebra $A$ which is $\mathbb{N}$-graded (namely $A =
\bigoplus_{i = 0}^{+\infty} A_i$) with $A_0 = K$, such that $L$ is
the fraction field of $A$. We can talk about homogeneous elements in
$L$ (i.e., quotients of two homogeneous elements of $A$, with
naturally defined degree). Now suppose $L$ is $K$-rational of
transcendence degree $n$, and can be generated by $n$ homogeneous
elements, say $y_1,\dots,y_n$, we may assume that all the $y_i$'s
have the same degree $e > 0$ (cf. \cite{Kem}, prop. 1.1 and its
proof). Then every homogeneous element $f \in L$ is of degree $ke$
for some $k \in \mathbb{Z}$, and moreover $f$ can be written
(uniquely) as a homogeneous rational function of $y_1,\dots,y_n$ of
degree $k$. In fact, $L = K(y_1,\dots,y_n) =
K(y_1,\frac{y_2}{y_1},\dots,\frac{y_n}{y_1})$ and the subfield $L_0
\subset L$ formed by homogeneous elements of degree 0 is just
$K(\frac{y_2}{y_1},\dots,\frac{y_n}{y_1})$. Now $\frac{f}{y_1^k} \in
L_0$ is a rational function of
$\frac{y_2}{y_1},\dots,\frac{y_n}{y_1}$, the result follows.

\begin{proposition}
Let $G$ be a group acting on the rational function field
$K(x_1,\dots,x_n)$ by some linear representation, namely, for any
$\sigma \in G$, we have $\sigma(x_i) = \sum_{j=1}^n a_{ij}(\sigma)
x_j$ with $a_{ij}(\sigma) \in K$. Let $x_{n+1}$ be another variable
with trivial $G$ action, and assume that $K(x_1,\dots,x_n)^G$ is
$K$-rational. Then $K(x_1,\dots,x_n,x_{n+1})^G$ is also
$K$-rational, and can be generated by $n+1$ homogeneous rational
functions of degree $1$.
\end{proposition}

To see this, just notice that $K(x_1,\dots,x_n,x_{n+1})^G =
K\big(\frac{x_1}{x_{n+1}},\dots,\frac{x_n}{x_{n+1}},x_{n+1}\big)^G
=$
$K\big(\frac{x_1}{x_{n+1}},\dots,\frac{x_n}{x_{n+1}}\big)^G(x_{n+1})$.
Now $K\big(\frac{x_1}{x_{n+1}},\dots,\frac{x_n}{x_{n+1}}\big)^G$ is
$K$-rational, and generated by $n$ elements $f_1,\dots,f_n$ which
are homogeneous of degree $0$, so $K(x_1,\dots,x_n,x_{n+1})^G$ is
generated by $x_{n+1},x_{n+1}f_1,\dots,x_{n+1}f_n$, all homogeneous
of degree $1$.

\bigskip

The following is a simple rationality criterion to be used later.

\begin{lemma}
Let $L$ be an extension field of $K$ generated by $x_1,\dots,x_n$
and $y$, with only one relation $$ y^m = f(x_1,\dots,x_n). $$
Assuming $f$ is a homogeneous rational function of $x_1,\dots,x_n$
of degree $k$, and $m,k$ are coprime. Then $L$ is $K$-rational.
\;{\rm (One can formulate a more general statement, but this form is
enough for our purpose.)}
\end{lemma}

The proof is easy: we let $z_1 = \frac{x_1}{x_n}, \dots, z_{n-1} =
\frac{x_{n-1}}{x_n}$, then $L = K(z_1,\dots,z_{n-1},$ $x_n,y)$ with
one relation $y^m x_n^{-k} = f(z_1,\dots,z_{n-1},1)$. Let $\xi = y^m
x_n^{-k}$ and $\eta = y^a x_n^b$ for some integers $a,b$ satisfying
$ka+mb=1$. It is then clear that $L = K(z_1,\dots,z_{n-1},x_n,y) =
K(z_1,\dots,z_{n-1},\xi,\eta) = K(z_1,\dots,z_{n-1},\eta)$.

\bigskip

In the case of positive characteristic, there is the following basic
result:

\begin{theorem} {\rm (\cite{Ku}, \cite{Mi})} Let $K$ be a field of
characteristic $p > 0$, and $G$ be a $p$-group acting on
$K(x_1,\dots,x_n)$ by some linear representation. Then the fixed
field $K(x_1,\dots,x_n)^G$ is $K$-rational.
\end{theorem}

For elementary 2-groups of order $4$ and $8$, a simple system of
generators for the fixed field is provided by:

\begin{proposition} Let $K$ be a field of characteristic 2.

{\rm (1)} Let $V_4 \subset S_4$ be the Klein 4 group generated by
$(1,2)(3,4),(1,3)(2,4)$, and act on $x_1,x_2,x_3,x_4$ by
permutations, then $K(x_1,x_2,x_3,x_4)^{V_4}$ is generated by $$ v_1
= x_1+x_2+x_3+x_4, \quad v_2 = x_1x_2+x_3x_4, \quad v_3 =
x_1x_3+x_2x_4, \quad v_4 = x_1x_4+x_2x_3. $$

{\rm (2)} Let $V_8 \subset S_8$ be the group generated by
$(1,2)(3,4)(5,6)(7,8),(1,3)(2,4)(5,7)(6,8)$, $(1,5)(2,6)(3,7)(4,8)$,
and act on $x_1,\dots,x_8$ by permutations, then
$K(x_1,\dots,x_8)^{V_8}$ is generated by $$
\begin{array}{l} w_1 = x_1+x_2+x_3+x_4+x_5+x_6+x_7+x_8, \quad
w_2 = x_1x_2+x_3x_4+x_5x_6+x_7x_8, \\
w_3 = x_1x_3+x_2x_4+x_5x_7+x_6x_8, \quad
w_4 = x_1x_4+x_2x_3+x_5x_8+x_6x_7, \\
w_5 = x_1x_5+x_2x_6+x_3x_7+x_4x_8, \quad
w_6 = x_1x_6+x_2x_5+x_3x_8+x_4x_7, \\
w_7 = x_1x_7+x_2x_8+x_3x_5+x_4x_6, \quad w_8 =
x_1x_8+x_2x_7+x_3x_6+x_4x_5.
\end{array} $$ {\rm (It is possible to formulate a similar result for any
elementary 2-group, but the above two cases are enough for our
purpose.)}
\end{proposition}

Proof: It is easily checked that these elements are fixed by
corresponding groups, so we are reduced to show that the degree of
extensions are 4 and 8 respectively.

Let $u_1 = x_1+x_2, u_2 = x_1x_2$, then $[K(x_1,x_2):K(u_1,u_2)] =
2$. For (1), we use the following identities:
\begin{align*} u_1^2 +
v_1 u_1 + v_3 + v_4 = 0 , \tag{i} \\ (v_3+v_4)^2 u_2 + u_1^4 u_2 +
v_2 u_1^4 + v_3 v_4 u_1^2 = 0 , \tag{ii}
\end{align*}
which show that
$$
\begin{array}{ll}
& [K(x_1,x_2,x_3,x_4) : K(v_1,v_2,v_3,v_4)] \\
=& [K(x_1,x_2,v_3,v_4) : K(v_1,v_2,v_3,v_4)] \quad \text{(since
$x_3 = \frac{x_1 v_3+x_2 v_4}{(x_1+x_2)^2},\, x_4 = \frac{x_2 v_3+x_1 v_4}{(x_1+x_2)^2}$)}\\
=& 2\,[K(u_1,u_2,v_3,v_4) : K(v_1,v_2,v_3,v_4)] \quad
\text{($v_1,v_2 \in K(u_1,u_2,v_3,v_4)$
follows from (i), (ii))} \\
=& 2\,[K(u_1,v_2,v_3,v_4) : K(v_1,v_2,v_3,v_4)] \quad \text{(by (ii) above)} \\
=& 4 \quad \text{(by (i) above).}
\end{array}
$$
For (2), we use the following identities:
\begin{align*} v_1^2 + w_1
v_1 + w_5 + w_6 + w_7 + w_8 = 0 , \tag{iii} \\
(w_5 + w_6 + w_7 + w_8)^2 v_2 + v_1^4 v_2 + w_2 v_1^4 + (w_5 w_6 +
w_7 w_8) v_1^2 = 0 , \tag{iv} \\
(w_5 + w_6 + w_7 + w_8)^2 v_3 + v_1^4 v_3 + w_3 v_1^4 + (w_5 w_7 +
w_6 w_8) v_1^2 = 0 , \tag{v} \\
(w_5 + w_6 + w_7 + w_8)^2 v_4 + v_1^4 v_4 + w_4 v_1^4 + (w_5 w_8 +
w_6 w_7) v_1^2 = 0 , \tag{vi}
\end{align*}
which show that $$
\begin{array}{ll}
& [K(x_1,x_2,x_3,x_4,x_5,x_6,x_7,x_8) :
K(w_1,w_2,w_3,w_4,w_5,w_6,w_7,w_8)] \\
=& [K(x_1,x_2,x_3,x_4,w_5,w_6,w_7,w_8) :
K(w_1,w_2,w_3,w_4,w_5,w_6,w_7,w_8)] \\
=& 4\,[K(v_1,v_2,v_3,v_4,w_5,w_6,w_7,w_8) :
K(w_1,w_2,w_3,w_4,w_5,w_6,w_7,w_8)] \quad \text{(by (1))} \\
=& 4\,[K(v_1,w_2,w_3,w_4,w_5,w_6,w_7,w_8) :
K(w_1,w_2,w_3,w_4,w_5,w_6,w_7,w_8)] \\ & \hfill
\text{(by (iv), (v), (vi) above)} \\
=& 8 \quad \text{(by (iii) above)}
\end{array}
$$
The proof is completed.

\bigskip

We will also make use of the following concrete result (valid for
any field):

\begin{theorem} {\rm (\cite{Mas}, Th. 3, \cite{HoK}, Th. 2.2)}
Let $K$ be a field, $x_1,x_2,x_3$ be variables, and $S_3$ acts on
$K(x_1,x_2,x_3)$ by permutation of variables. For the alternating
subgroup $A_3 = \langle (1,2,3)\rangle$, the fixed field
$K(x_1,x_2,x_3)^{A_3}$ is $K$-rational, and generated by $$
x_1+x_2+x_3, \quad
\frac{x_1x_2^2+x_2x_3^2+x_3x_1^2-3x_1x_2x_3}{x_1^2+x_2^2+x_3^2-x_1x_2-x_2x_3-x_3x_1},
\quad
\frac{x_1x_3^2+x_2x_1^2+x_3x_2^2-3x_1x_2x_3}{x_1^2+x_2^2+x_3^2-x_1x_2-x_2x_3-x_3x_1}.
$$
\end{theorem}

Finally, in counting the extension degree, the following lemma is
frequently used.

\begin{lemma} Let $K(x_1,\dots,x_n)$ be a rational function field on
variables $x_1,\dots,x_n$, and $$ f_i = \prod_{k=1} x_k^{a_{ik}}
\quad (1 \le i \le n) $$ a set of monomials. Assume the determinant
$D$ of the matrix $(a_{ik})$ is nonzero, then $K(x_1,\dots,x_n)$ is
a finite extension of $K(f_1,\dots,f_n)$, with extension degree $=
|D|$.
\end{lemma}

The proof is easy, omitted.


\section{Subgroups of $S_8$}

\cite{CHM} contains a complete list of all transitive subgroups of
$S_8$ not equal to $S_8$ and $A_8$ up to conjugation. For the later
use, we first introduce some notations:
$$ \begin{array}{llllll}
\sigma_1 &=& (1,2)(3,4)(5,6)(7,8), & \qquad \sigma_2 &=&
(1,3)(2,4)(5,7)(6,8),
\\ \kappa &=& (1,5)(2,6)(3,7)(4,8), & \qquad \widetilde\kappa &=&
(1,5)(2,6)(3,8)(4,7),
\\ \kappa' &=& (1,5)(2,6)(3,8,4,7), & \qquad \kappa^\circ &=&
(1,6)(2,5)(3,7)(4,8)
\\ \varphi_1 &=& (2,3,4)(7,6,5), & \qquad \widetilde\varphi_1 &=&
(2,3,4)(6,7,8),
\\ \varphi_2 &=& (2,3)(7,6), & \qquad \widetilde\varphi_2 &=& (3,4)(5,6),
\\ \psi_1 &=& (2,7)(3,4,6,5), & \qquad \psi_2 &=& (2,4)(5,6,7,8),
\\ \theta &=& (2,3,7,4,5,6,8), & \qquad \widetilde\theta &=& (1,2,3,4,5,6,7),
\\ \Theta &=& (1,2,3,4,5,6,7,8), & \qquad \Phi &=& (1,2,4)(5,6,8),
\\ \Psi &=& (1,2,3,4)(5,6,7,8), & \qquad \widetilde\Psi &=& (1,3,5,7)(2,4,6,8).
\end{array} $$
We list the groups according to their orders as follows:
\begin{itemize}
\item Groups of order $8$:
$$ \begin{array}{lll}
G_1 &=& \langle \Theta\rangle \\
G_2 &=& \langle \Psi,\; \kappa\rangle \\
G_3 &=& \langle \sigma_1,\; \sigma_2,\; \kappa\rangle \\
G_4 &=& \langle \Psi,\; (1,8)(2,7)(3,6)(4,5)\rangle \\
G_5 &=& \langle \Psi,\; (1,5,3,7)(2,8,4,6)\rangle
\end{array} $$
\item Groups of order $16$:
$$ \begin{array}{lll}
G_6 &=& \langle \Theta,\; (1,8)(2,7)(3,6)(4,5)\rangle \\
G_7 &=& \langle \Theta,\; (2,6)(4,8)\rangle \\
G_8 &=& \langle \Theta,\; (2,4)(3,7)(6,8)\rangle \\
G_9 &=& \langle \sigma_1,\; \sigma_2,\; \kappa,\; (5,6)(7,8)\rangle \\
G_{10} &=& \langle (2,6)(4,8),\; \Psi\rangle \\
G_{11} &=& \langle (2,6)(4,8),\; \widetilde\Psi,\;
(1,2,5,6)(3,4,7,8)\rangle
\end{array} $$
\item Groups of order $24$:
$$ \begin{array}{lll}
G_{12} &=& \langle \Phi,\; \widetilde\Psi\rangle \\
G_{13} &=& \langle \sigma_1,\; \sigma_2,\; \kappa,\; \varphi_1\rangle \\
G_{14} &=& \langle \sigma_2,\; \widetilde\varphi_1,\;
\kappa^\circ\rangle
\end{array} $$
\item Groups of order $32$:
$$ \begin{array}{lll}
G_{15} &=& \langle \Theta,\; (2,6)(4,8),\; (1,8)(2,7)(3,6)(4,5)\rangle \\
G_{16} &=& \langle \Theta,\; (3,7)(4,8)\rangle \\
G_{17} &=& \langle (1,2,3,4),\; \kappa\rangle \\
G_{18} &=& \langle \sigma_1,\; \sigma_2,\; \kappa,\; (5,6)(7,8),\; (5,7)(6,8)\rangle \\
G_{19} &=& \langle \sigma_1,\; \sigma_2,\; \kappa,\; \psi_2\rangle \\
G_{20} &=& \langle (3,7)(4,8),\; \Psi\rangle \\
G_{21} &=& \langle (2,6)(4,8),\; (1,2,5,6)(3,4)(7,8),\; \sigma_2\rangle \\
G_{22} &=& \langle \sigma_1,\; \sigma_2,\; \kappa,\;
\widetilde\varphi_2,\; (3,4)(7,8)\rangle
\end{array} $$
\item Groups of order $48$:
$$ \begin{array}{lll}
G_{23} &=& \langle \Theta,\; \Phi\rangle \\
G_{24} &=& \langle \sigma_1,\; \sigma_2,\; \kappa,\; \varphi_1,\;
\widetilde\varphi_2\rangle
\end{array} $$
\item Groups of order $56$:
$$ \begin{array}{lll}
G_{25} &=& \langle \sigma_1,\; \sigma_2,\; \kappa,\; \theta\rangle
\end{array} $$
\item Groups of order $64$:
$$ \begin{array}{lll}
G_{26} &=& \langle \Theta,\; (1,5)(2,6),\; (1,5)(2,8)(4,6)\rangle \\
G_{27} &=& \langle (1,5),\; \Psi\rangle \\
G_{28} &=& \langle (3,7)(4,8),\; (2,4)(6,8),\; \Theta\rangle \\
G_{29} &=& \langle \sigma_1,\; \sigma_2,\; \kappa,\; \psi_2,\; (2,4)(6,8)\rangle \\
G_{30} &=& \langle (3,7)(4,8),\; (1,5)(2,4)(6,8),\; \Psi\rangle \\
G_{31} &=& \langle (1,5),\; \sigma_1,\; \sigma_2\rangle
\end{array} $$
\item Groups of order $96$:
$$ \begin{array}{lll}
G_{32} &=& \langle \sigma_1,\; \sigma_2,\; \kappa,\; \varphi_1,\; \psi_1^2\rangle \\
G_{33} &=& \langle \sigma_1,\; \sigma_2,\; \kappa,\; \varphi_1,\; (5,7)(6,8)\rangle \\
G_{34} &=& \langle (1,2)(3,4),\; \widetilde\varphi_1,\;
\widetilde\kappa\rangle
\end{array} $$
\item Groups of order $128$:
$$ \begin{array}{lll}
G_{35} &=& \langle (1,5),\; (2,4)(6,8),\; \Psi\rangle
\end{array} $$
\item Groups of order $168$:
$$ \begin{array}{lll}
G_{36} &=& \langle \sigma_1,\; \sigma_2,\; \kappa,\; \theta,\; \varphi_1\rangle \\
G_{37} &=& \langle \widetilde\theta,\; (2,3,5)(4,7,6),\;
(1,8)(2,7)(3,4)(5,6)\rangle
\end{array} $$
\item Groups of order $192$:
$$ \begin{array}{lll}
G_{38} &=& \langle (1,5),\; \sigma_1,\; \widetilde\varphi_1\rangle \\
G_{39} &=& \langle \sigma_1,\; \sigma_2,\; \kappa,\; \varphi_1,\; \psi_1\rangle \\
G_{40} &=& \langle (1,5)(2,6),\; \sigma_1,\; \widetilde\varphi_1,\; (1,5)(3,4)(7,8)\rangle \\
G_{41} &=& \langle \sigma_1,\; \sigma_2,\; \kappa,\; \varphi_1,\;
\psi_2\rangle
\end{array} $$
\item Groups of order $>192$:
$$ \begin{array}{llll}
G_{42} &=& \langle (1,3)(2,4),\; (2,3,4),\; \kappa\rangle & \text{of order}\; 288 \\
G_{43} &=& \langle \widetilde\theta,\; (1,8)(2,7)(3,4)(5,6),\; (2,4,3,7,5,6)\rangle & \text{of order}\; 336 \\
G_{44} &=& \langle (1,5),\; (1,2)(5,6),\; \Psi\rangle & \text{of order}\; 384 \\
G_{45} &=& \langle (1,3)(2,4),\; (2,3,4),\; (1,2)(5,6),\; \kappa\rangle & \text{of order}\; 576 \\
G_{46} &=& \langle (1,3)(2,4),\; (2,3,4),\; (1,2)(5,6),\; \kappa'\rangle & \text{of order}\; 576 \\
G_{47} &=& \langle (1,2,3,4),\; (3,4),\; \kappa\rangle & \text{of order}\; 1152 \\
G_{48} &=& \langle \sigma_1,\; \sigma_2,\; \theta,\; \kappa,\;
\varphi_1,\; \varphi_2\rangle & \text{of order}\; 1344
\end{array} $$
\end{itemize}

Among these groups, only $G_{37}, G_{43}$ and $G_{48}$ are
unsolvable.

\bigskip

In the following, we will divide these groups into 4 types (exclude
one group $G_{43}$, which we don't know the answer), based on the
methods employed for the rationality proof.

Type A contains 23 groups, numbered by 1$\sim$12, 17, 18, 23, 27,
31, 35, 37, 38, 42, 44, 47. For these groups, the results are well
known, and distributed in literatures.

Type B contains 16 groups, numbered by 13$\sim$16, 19$\sim$22, 24,
26, 28$\sim$30, 32, 39, 40. For these groups, we can reduce the
rationality problem to corresponding problems for 3-dimensional
monomial actions, then use the methods of \cite{HoKiYa}. Some
efforts is needed here for the characteristic 2 cases.

Type C contains 5 groups, numbered by 33, 34, 41, 45, 46. For these
groups, the rationality problem is reduced to corresponding problems
for subgroups of $S_6$, then use Theorem 1.1 (2).

Type D contains 3 groups, numbered by 25, 36, 48. For these groups,
the rationality problem is reduced to corresponding problems for
subgroups of $S_7$, then use Theorem 1.1 (3).

\section{Various known cases}

We collect them into 3 classes:
\begin{itemize}
\item Wreath products, totally of 9 groups:
$G_{17} = C_4 \wr C_2, \; G_{18} = V_4 \wr C_2, \; G_{42} = A_4 \wr
C_2, \; G_{47} = S_4 \wr C_2, \; G_{27} = C_2 \wr C_4, \; G_{31} =
C_2 \wr V_4, \; G_{35} = C_2 \wr D_4, \; G_{38} = C_2 \wr A_4, \;
G_{44} = C_2 \wr S_4$. 
\item Groups of order $8$ and $16$, totally of 11 groups:
$G_1,G_2,G_3,G_4,G_5$ (of order $8$) and
$G_6,G_7,G_8,G_9,G_{10},G_{11}$ (of order $16$).
\item Finite linear groups: $G_{12} \simeq {\rm SL}(2,3)$, $G_{23}
\simeq {\rm GL}(2,3)$ and $G_{37} \simeq {\rm PSL}(2,7)$ (note that
$G_{43} \simeq {\rm PGL}(2,7)$ also belongs to this class, and it is
the only case remains unknown).
\end{itemize}

The first class is taken care by Theorem 2.3 and Theorem 1.1 (1).

\bigskip

In the second class, the results are essentially contained in
\cite{Len}, \cite{ChHuK} and \cite{Ka}. Actually, \cite{Len}
contains a complete treatment of the rationality problem for $K(x_g
: g\in G)^G$ (known as Noether's problem) when $G$ is an abelian
group, so the cases $G_1 \simeq C_8,\, G_2 \simeq C_2 \times C_4,\,
G_3 \simeq C_2 \times C_2 \times C_2$ follow directly from the Main
Theorem of \cite{Len}. Also, \cite{ChHuK} contains results of $K(x_g
: g\in G)^G$ for nonabelian group $G$ of order $8$, see
\cite{ChHuK}, Prop. 2.6 for the case $G_4 \simeq D_4$ and Th. 2.7
for the case $G_5 \simeq Q_8$.

On the other hand, the results of \cite{ChHuK} and \cite{Ka} on
Noether's problem (i.e. rationality problem of $K(x_g : g\in G)^G$)
for nonabelian group $G$ of order $16$ cannot apply directly to our
problem, but the methods of proofs used there essentially give
solutions for our cases. In fact,

$G_6 \simeq D_8$ is treated in \cite{ChHuK}, Th. 3.1 (by taking
$\sigma = \Theta, \tau = (1,8)(2,7)(3,6)(4,5)\Theta$, and replacing
$8$ by $0$);

$G_7 \simeq$ modular group of order $16$ is treated in \cite{ChHuK},
Th. 3.3 (by taking $\sigma = \Theta, \tau =
\Theta^{-1}(2,6)(4,8)\Theta$, and replacing $8$ by $0$);

$G_8 \simeq$ quasi-dihedral group of order $16$ is treated in
\cite{ChHuK}, Th. 3.2 (by taking $\sigma = \Theta, \tau =
\Theta^{-1}(2,4)(3,7)(6,8)\Theta$, and replacing $8$ by $0$);

$G_9 \simeq C_2 \times D_4$ is the group (VI) in \cite{Ka}, p. 305
and treated in p. 307-308 (by taking $\sigma = \kappa (5,6)(7,8) =
(1,5,2,6)(3,7,4,8),\, \tau = (5,6)(7,8),\, \lambda = \sigma_2$, and
letting
$$
\begin{array}{llll} x_1 = x_e + x_\tau, & x_2 = x_{\sigma^2} +
x_{\sigma^2 \tau}, & x_3 = x_\lambda + x_{\lambda \tau}, & x_4 =
x_{\lambda \sigma^2} + x_{\lambda \sigma^2 \tau}, \\ x_5 = x_\sigma
+ x_{\sigma \tau}, & x_6 = x_{\sigma^3} + x_{\sigma^3 \tau}, & x_7 =
x_{\lambda \sigma} + x_{\lambda \sigma \tau}, & x_8 = x_{\lambda
\sigma^3} + x_{\lambda \sigma^3 \tau}, \\ X_1 = x_1 - x_2, & X_2 =
x_5 - x_6, & X_3 = x_3 - x_4, & X_4 = x_7 - x_8,
\end{array}
$$ these $X_1, X_2, X_3, X_4$ are used in the proof there);

$G_{10}$ is the group (IX) in \cite{Ka}, p. 305 and treated in p.
310-311 (by taking $\sigma = \Psi,\, \tau = \Psi^{-1} (2,6)(4,8)
\Psi (2,6)(4,8) = \kappa,\, \lambda = \Psi^2 (2,6)(4,8) =
(1,3)(5,7)(2,8)(4,6)$, and letting
$$
\begin{array}{llll}
x_1 = x_e + x_{\sigma^2 \lambda}, & x_2 = x_\sigma + x_{\sigma^3
\lambda}, & x_3 = x_{\sigma^2} + x_\lambda, & x_4 = x_{\sigma^3} +
x_{\sigma \lambda}, \\ x_5 = x_\tau + x_{\tau \sigma^2 \lambda}, &
x_6 = x_{\tau \sigma} + x_{\tau \sigma^3 \lambda}, & x_7 = x_{\tau
\sigma^2} + x_{\tau \lambda}, & x_8 = x_{\tau \sigma^3} + x_{\tau
\sigma \lambda}, \\ X_1 = x_1 - x_3, & X_2 = x_2 - x_4, & X_3 = x_5
- x_7, & X_4 = x_6 - x_8,
\end{array}
$$ these $X_1, X_2, X_3, X_4$ are used in the proof there);

$G_{11}$ is the group (V) in \cite{Ka}, p. 305 and treated in p.
305-307 (by taking $\sigma = \widetilde\Psi,\, \tau = (2,6)(4,8),\,
\lambda = (1,2,5,6)(3,4,7,8) \widetilde\Psi = (1,4)(3,6)(5,8)(7,2)$,
and letting
$$
\begin{array}{llll}
x_1 = x_e + x_\tau, & x_2 = x_{\lambda \sigma^3} + x_{\lambda
\sigma^3 \tau}, & x_3 = x_\sigma + x_{\sigma \tau}, & x_4 =
x_\lambda + x_{\lambda \tau},
\\ x_5 = x_{\sigma^2} + x_{\sigma^2 \tau}, & x_6 = x_{\lambda \sigma}
+ x_{\lambda \sigma \tau}, & x_7 = x_{\sigma^3} + x_{\sigma^3 \tau}, &
x_8 = x_{\lambda \sigma^2} + x_{\lambda \sigma^2 \tau}, \\
X_1 = x_1 - x_5, & X_2 = x_3 - x_7, & X_3 = x_4 - x_8, & X_4 = x_6 -
x_2,
\end{array}
$$ these $X_1, X_2, X_3, X_4$ are used in the proof there).

\bigskip

In the third class, $G_{37} \simeq {\rm PSL}(2,7)$ can be treated
exactly the same way as in Theorem 1.1 (3), detail can be found in
\cite{KW}, proof of Th. 4.7, and we avoid the repetition here.

The first two groups $G_{12}$ and $G_{23}$ are treated by \cite{Ri}
and Plans \cite{Pl}, we reformulate their results here to fit into
our framework:

$\mathbb{F}_3^2$ has $8$ nonzero elements, and the natural action of
${\rm GL}(2,3)$ on $\mathbb{F}_3^2 \smallsetminus
\big\{\binom00\big\}$ corresponds to the action of $G_{23}$ on
$x_i$'s under the following correspondence
$$
\begin{pmatrix} 1 \\ 0 \end{pmatrix} \mapsto x_1, \;\;
\begin{pmatrix} 1 \\ -1 \end{pmatrix} \mapsto x_2, \;\;
\begin{pmatrix} 0 \\ 1 \end{pmatrix} \mapsto x_3, \;\;
\begin{pmatrix} 1 \\ 1 \end{pmatrix} \mapsto x_4, $$ $$
\begin{pmatrix} -1 \\ 0 \end{pmatrix} \mapsto x_5, \;\;
\begin{pmatrix} -1 \\ 1 \end{pmatrix} \mapsto x_6, \;\;
\begin{pmatrix} 0 \\ -1 \end{pmatrix} \mapsto x_7, \;\;
\begin{pmatrix} -1 \\ -1 \end{pmatrix} \mapsto x_8, $$
and the subgroup ${\rm SL}(2,3)$ corresponds to $G_{12}$. More
precisely, one checks that
$$
\begin{pmatrix} 1 & 1 \\ -1 & 1 \end{pmatrix}
\longleftrightarrow \Theta, \quad
\begin{pmatrix} 1 & 0 \\ -1 & 1 \end{pmatrix}
\longleftrightarrow \Phi, \quad
\begin{pmatrix} 0 & -1 \\ 1 & 0 \end{pmatrix}
\longleftrightarrow \widetilde\Psi.
$$

As the methods are the same, we will consider $G_{23}$ only. When
${\rm char}(K) \ne 2$, by changing variables
$$ \begin{array}{llll}
y_1 = x_1-x_5, & y_2 = x_2-x_6, & y_3 = x_3-x_7, & y_4 = x_4-x_8, \\
y_5 = x_1+x_5, & y_6 = x_2+x_6, & y_7 = x_3+x_7, & y_8 = x_4+x_8,
\end{array} $$
we get $$ \begin{array}{lll} \Theta &:& y_1 \mapsto y_2 \mapsto y_3
\mapsto y_4 \mapsto -y_1,\; y_5 \mapsto y_6 \mapsto y_7 \mapsto y_8
\mapsto y_5, \\ \Phi &:& y_1 \mapsto y_2 \mapsto y_4 \mapsto y_1,\;
y_5 \mapsto y_6 \mapsto y_8 \mapsto y_5,\; \text{and}\; y_3, y_7
\;\text{fixed}.
\end{array} $$
This shows that $K(y_1,y_2,y_3,y_4)$ is stable under the action of
$G_{23}$, and the restricted action of $G_{23}$ is stable by trivial
reason (if $g(x_i-x_j) = x_i-x_j$, we must have $g(x_i)=x_i$ and
$g(x_j)=x_j$, since ${\rm char}(K) \ne 2$). So by Theorem 2.1 (1),
we are reduced to prove the rationality of
$K(y_1,y_2,y_3,y_4)^{G_{23}}$ and $K(y_1,y_2,y_3,y_4)^{G_{12}}$.

Now $G_{23}$ (and $G_{12}$) has a normal subgroup $Q = \langle
\widetilde\Psi,\; \Phi^{-1}\widetilde\Psi\Phi\rangle$ (note that
$\Phi^{-1}\widetilde\Psi\Phi = (1,2,5,6)(4,3,8,7)$), which is
isomorphic to the quaternion group $Q_8$, and we see that
$$ \begin{array}{lll} \widetilde\Psi &:& y_1
\mapsto y_3 \mapsto -y_1, \quad y_2
\mapsto y_4 \mapsto -y_2, \\
\Phi^{-1}\widetilde\Psi\Phi &:& y_1 \mapsto y_2 \mapsto -y_1, \quad
y_3 \mapsto -y_4 \mapsto -y_3.
\end{array} $$
There is a ``canonical" basis for the fixed field
$K(y_1,y_2,y_3,y_4)^Q$ (cf. \cite{Pl}): $$ z_1 =
\frac{y_1y_2-y_3y_4}{y_2y_4+y_1y_3}, \quad z_2 =
\frac{y_2y_4-y_1y_3}{y_4y_1+y_2y_3}, \quad z_3 =
\frac{y_4y_1-y_2y_3}{y_1y_2+y_3y_4}, $$
$$ z_4 = y_1^2+y_2^2+y_3^2+y_4^2. $$ 
And we find $$ \begin{array}{lll} \Theta &:& \displaystyle z_1
\mapsto \frac1{z_2}, \quad z_2 \mapsto \frac1{z_1}, \quad z_3
\mapsto \frac1{z_3}, \quad z_4 \mapsto z_4, \\ \Phi &:& z_1 \mapsto
z_2 \mapsto z_3 \mapsto z_1, \quad z_4 \mapsto z_4.
\end{array} $$
Now we only need to apply Theorem 2.4 (2).

When ${\rm char}(K) = 2$, there is a similar argument. Indeed, we
may use Proposition 2.2 to get the fixed field of $\kappa$, namely
$K(x_1,\dots,x_8)^{\langle\kappa\rangle} = K(y_1,\dots,y_8)$ with
$$
\begin{array}{llll} y_1 = x_1+x_5, & y_2 = x_2+x_6, & y_3 = x_3+x_7,
& y_4 = x_4+x_8, \\ y_5 = x_1x_5, & y_6 = x_1x_6+x_5x_2, & y_7 =
x_1x_7+x_5x_3, & y_8 = x_1x_8+x_5x_4.
\end{array} $$
Then the action of $Q$ on $K(y_1,\dots,y_8)$ is reduced to a $V_4$
action as follows:
$$ \begin{array}{lll} \widetilde\Psi &:& \displaystyle y_1
\leftrightarrow y_3, \quad y_2 \leftrightarrow y_4, \quad y_5
\mapsto \frac{y_3^2}{y_1^2} y_5 + \frac{y_3}{y_1} y_7 +
\frac{1}{y_1^2} y_7^2, \\
&& \displaystyle y_6 \mapsto \frac{y_4}{y_1} y_7 + \frac{y_3}{y_1}
y_8, \quad y_7 \mapsto y_1y_3 + y_7, \quad y_8 \mapsto y_2y_3 +
\frac{y_3}{y_1} y_6 + \frac{y_2}{y_1} y_7, \\
\Phi^{-1}\widetilde\Psi\Phi &:& \displaystyle y_1 \leftrightarrow
y_2, \quad y_3 \leftrightarrow y_4, \quad y_5 \mapsto
\frac{y_2^2}{y_1^2} y_5 + \frac{y_2}{y_1} y_6 + \frac{1}{y_1^2}
y_6^2, \\ && \displaystyle y_6 \mapsto y_1y_2 + y_6, \quad y_7
\mapsto y_2y_4 + \frac{y_4}{y_1} y_6 + \frac{y_2}{y_1} y_8, \quad
y_8 \mapsto \frac{y_3}{y_1} y_6 + \frac{y_2}{y_1} y_7.
\end{array} $$
Replace $y_5,y_6,y_7,y_8$ by some elements which are invariant under
$Q$ as follows:
$$ \begin{array}{l}
\displaystyle z_5 = \frac{y_1^2+y_2^2+y_3^2+y_4^2}{y_1^2}y_5+
\frac{y_2}{y_1}y_6+\frac{1}{y_1^2}y_6^2+
\frac{y_3}{y_1}y_7+\frac{1}{y_1^2}y_7^2+
\frac{y_4}{y_1}y_8+\frac{1}{y_1^2}y_8^2, \\
\displaystyle z_6 = y_1^2 y_2+y_3^2
y_4+(y_1+y_2)y_6+\frac{y_4^2+y_3y_4}{y_1}y_7+
\frac{y_3^2+y_3y_4}{y_1}y_8, \\
\displaystyle z_7 = y_1^2 y_3+y_4^2
y_2+\frac{y_4^2+y_2y_4}{y_1}y_6+(y_1+y_3)y_7+
\frac{y_2^2+y_2y_4}{y_1}y_8, \\
\displaystyle z_8 = y_1^2 y_4+y_2^2 y_3+\frac{y_3^2+y_2y_3}{y_1}y_6+
\frac{y_2^2+y_2y_3}{y_1}y_7+(y_1+y_4)y_8,
\end{array} $$
we see immediately that $K(y_1,y_2,y_3,y_4)(y_5,y_6,y_7,y_8) =
K(y_1,y_2,y_3,y_4)(z_5,z_6,z_7,z_8)$. So by Proposition 2.9 (1), we
have
\begin{multline*}
K(y_1,y_2,y_3,y_4,y_5,y_6,y_7,y_8)^Q =
K(y_1,y_2,y_3,y_4)(z_5,z_6,z_7,z_8)^Q \\ =
K(y_1,y_2,y_3,y_4)^Q(z_5,z_6,z_7,z_8) =
K(z_1,z_2,z_3,z_4,z_5,z_6,z_7,z_8) \end{multline*} with $$ z_1 =
y_1+y_2+y_3+y_4, \quad z_2 = y_1 y_2+y_3 y_4, \quad z_3 = y_1
y_3+y_2 y_4, \quad z_4 = y_2 y_3+y_1 y_4. $$ Furthermore, the
actions of $\Theta$ and $\Phi$ on $z_i$'s are given by $$
\begin{array}{lll} \Theta &:& z_2 \leftrightarrow z_4,
\quad z_6 \leftrightarrow z_8, \quad z_7 \mapsto z_7 + A_1, \;
\text{and} \; z_1, z_3, z_5 \; \text{fixed},
\\ \Phi &:& z_2 \mapsto z_3 \mapsto z_4 \mapsto z_2, \quad
z_6 \mapsto z_7 + A_1, \quad z_7 \mapsto z_8 \mapsto z_6 + A_2, \;
\text{and} \; z_1, z_5 \; \text{fixed},
\end{array} $$ where $A_1, A_2 \in K(z_1,z_2,z_3,z_4)$ \footnote{In fact,
$A_1 = z_1z_3+\frac{z_2z_3+z_3z_4+z_4z_2}{z_1}$ and $A_2 =
z_1z_2+\frac{z_2z_3+z_3z_4+z_4z_2}{z_1}$.}.

Now the action of $G_{23}/Q$ on $K(z_1,z_2,z_3,z_4)$ is faithful
(because $|G_{23}/Q| = 6$), and its action on $z_5,z_6,z_7,z_8$ is
affine-linear, with coefficients in $K(z_1,z_2,z_3,z_4)$, so by
Theorem 2.1 (1), $$ K(z_1,z_2,z_3,z_4,z_5,z_6,z_7,z_8)^{G_{23}/Q} =
K(z_1,z_2,z_3,z_4)^{G_{23}/Q}(w_5,w_6,w_7,w_8)
$$ for some $w_5,w_6,w_7,w_8$. From the above we see easily that
the action of $G_{23}/Q$ on $z_2,z_3,z_4$ is just the $S_3$ action,
so $K(z_1,z_2,z_3,z_4)^{G_{23}/Q} = K(w_1,w_2,w_3,w_4)$ with $w_1 =
z_1,\, w_2 = z_2+z_3+z_4,\, w_3 = z_2z_3+z_3z_4+z_4z_2,\, w_4 =
z_2z_3z_4$. (For $G_{12}$, the final step is an $A_3$ action, we may
use Theorem 2.10).

\begin{remark}
Another approach is by taking the following change of variables
$$ y_1 = \frac{x_1}{x_1+x_5}, \quad y_2 = \frac{x_2}{x_2+x_6}, \quad
y_3 = \frac{x_3}{x_3+x_7}, \quad y_4 = \frac{x_4}{x_4+x_8}, $$ $$
y_5 = x_1+x_5, \quad y_6 = x_2+x_6, \quad y_7 = x_3+x_7, \quad y_8 =
x_4+x_8,
$$ see \cite{KWZ}, section 3 for indications of this method.
\end{remark}

\section{Subgroups related to three-dimensional monomial
problems}

There are standard techniques to reduce the rationality problems to
lower dimensional ones, namely by finding some faithful
subrepresentation and using Theorem 1.1. This applies to our
problems on transitive subgroups of $S_8$ in many cases.

When ${\rm char}(K) \ne 2$, such groups are gathered into 4 sets as
follows:

\begin{itemize}
\item The first set contains 9 groups: $G_{14}, G_{15}, G_{16},
G_{20}, G_{21}, G_{26}, G_{28}, G_{30}, G_{40}$. For these groups,
we make variable change
$$ \begin{array}{llll}
y_1 = x_1-x_5, & y_2 = x_2-x_6, & y_3 = x_3-x_7, & y_4 = x_4-x_8, \\
y_5 = x_1+x_5, & y_6 = x_2+x_6, & y_7 = x_3+x_7, & y_8 = x_4+x_8.
\end{array} $$
\item The second set contains 4 groups: $G_{13}, G_{24}, G_{32},
G_{39}$. For these, we make variable change
$$ \begin{array}{llll}
y_1 = x_1-x_8, & y_2 = x_2-x_7, & y_3 = x_3-x_6, & y_4 = x_4-x_5, \\
y_5 = x_1+x_8, & y_6 = x_2+x_7, & y_7 = x_3+x_6, & y_8 = x_4+x_5.
\end{array} $$
\item The third set contains 2 groups: $G_{19}, G_{29}$. For these,
we make variable change
$$ \begin{array}{llll}
y_1 = x_1-x_3, & y_2 = x_2-x_4, & y_3 = x_5-x_7, & y_4 = x_6-x_8, \\
y_5 = x_1+x_3, & y_6 = x_2+x_4, & y_7 = x_5+x_7, & y_8 = x_6+x_8,
\end{array} $$
\item The fourth set contains 1 group: $G_{22}$.
For this, we make variable change
$$ \begin{array}{llll}
y_1 = x_1-x_2, & y_2 = x_3-x_4, & y_3 = x_5-x_6, & y_4 = x_7-x_8, \\
y_5 = x_1+x_2, & y_6 = x_3+x_4, & y_7 = x_5+x_6, & y_8 = x_7+x_8,
\end{array} $$
\end{itemize}
Then in each of these cases, one checks directly that the subfield
$K(y_1,y_2,y_3,y_4)$ is stable under the group action, and the
restricted action is faithful. So we are reduced to prove the
rationality of $K(y_1,y_2,y_3,y_4)^{G_i}$.

Letting in each case $$ z_1 = \frac{y_1}{y_4}, \quad z_2 =
\frac{y_2}{y_4}, \quad z_3 = \frac{y_3}{y_4}, $$ again by direct
computations, we see that $K(z_1,z_2,z_3)$ is stable under $G_i$,
and the action of $G_i$ on $K(z_1,z_2,z_3)$ is a monomial action.
And by Theorem 1.1 (2), we are reduced to consider
$K(z_1,z_2,z_3)^{G_i}$.

In the following, we will first use Theorem 2.4 and 2.5 to give the
rationality proofs under the assumption ${\rm char}(K) \ne 2$, our
discussion will be in several parts.

\bigskip

1. There are 3 groups whose actions on $K(z_1,z_2,z_3)$ are already
purely monomial, namely
$$
\begin{array}{ll}
G_{14} : & \left\{\begin{array}{ll} \displaystyle \sigma_2 \;:\; z_1
\mapsto \frac{z_3}{z_2}, \; z_2 \mapsto \frac{1}{z_2}, \; z_3
\mapsto \frac{z_1}{z_2}; & \quad \displaystyle \widetilde\varphi_1
\;:\; z_1 \mapsto \frac{z_1}{z_2}, \; z_2 \mapsto \frac{z_3}{z_2},
\; z_3 \mapsto \frac{1}{z_2};  \\ \displaystyle \kappa^\circ \;:\;
z_1 \mapsto z_2, \; z_2 \mapsto z_1, \; z_3 \mapsto z_2.
\end{array} \right. \\ \\
G_{13} : & \left\{\begin{array}{ll} \displaystyle \sigma_1 \;:\; z_1
\mapsto \frac{z_2}{z_3}, \; z_2 \mapsto \frac{z_1}{z_3}, \; z_3
\mapsto \frac{1}{z_3}; & \quad \displaystyle \sigma_2 \;:\; z_1
\mapsto \frac{z_3}{z_2}, \; z_2 \mapsto \frac{1}{z_2}, \; z_3
\mapsto \frac{z_1}{z_2}; \\
\displaystyle \kappa \;:\; z_1 \mapsto
\frac{1}{z_1}, \; z_2 \mapsto \frac{z_3}{z_1}, \; z_3 \mapsto
\frac{z_2}{z_1}; & \quad \displaystyle \varphi_1 \;:\; z_1 \mapsto
\frac{z_1}{z_2}, \; z_2 \mapsto \frac{z_3}{z_2}, \; z_3 \mapsto
\frac{1}{z_2}.
\end{array} \right. \\ \\
G_{24} : & \begin{array}{l} \sigma_1,\; \sigma_2,\; \kappa,\;
\varphi_1 \;\; \text{as in} \;\; G_{13}, \;\; \text{and} \;\;
\displaystyle \widetilde\varphi_2 \;:\; z_1 \mapsto \frac{z_1}{z_3},
\; z_2 \mapsto \frac{z_2}{z_3}, \; z_3 \mapsto \frac{1}{z_3}.
\end{array}
\end{array}
$$
So the results follows directly from Theorem 2.4 (2).

For other groups, we need to put the actions into reduced form, and
there is a uniform way to do this for most groups. In fact, one
shows easily that
\begin{itemize}
\item $\Lambda_1 = A_8 \cap \langle (1,5), (2,6), (3,7), (4,8)\rangle$
is a normal subgroup of order 8 of $G_{16}, G_{20},$ $G_{21},
G_{26}, G_{28}, G_{30}, G_{40}$,
\item $\Lambda_2 = A_8 \cap \langle (1,8), (2,7), (3,6), (4,5)\rangle$
is a normal subgroup of order 8 of $G_{32}, G_{39}$,
\item $\Lambda_3 = A_8 \cap \langle (1,3), (2,4), (5,7), (6,8)\rangle$
is a normal subgroup of order 8 of $G_{29}$,
\item $\Lambda_4 = A_8 \cap \langle (1,2), (3,4), (5,6), (7,8)\rangle$
is a normal subgroup of order 8 of $G_{29}$.
\end{itemize}
(The remaining 2 groups $G_{15}$ and $G_{19}$ will be discussed
later.) For each of these groups, the fixed field of $\Lambda_k$ can
be written as $K(z_1,z_2,z_3)^{\Lambda_k} = K(Z_1,Z_2,Z_3)$ with
$$ Z_1 = \frac{z_1z_2}{z_3}, \quad Z_2 = \frac{z_2z_3}{z_1}, \quad
Z_3 = \frac{z_3z_1}{z_2},
$$
and the action of $\overline{G}_i = G_i/(\text{some}\,\Lambda_k)$ on
$K(Z_1,Z_2,Z_3)$ is again monomial.

\bigskip

2. There are 5 groups $\overline{G}_i$ for which the actions on
$K(Z_1,Z_2,Z_3)$ are already purely monomial, namely
$$
\begin{array}{ll}
\overline{G}_{20} : & \displaystyle (3,7)(4,8) \;=\; {\rm Id}\,;
\quad \Psi \;:\; Z_1 \mapsto Z_2, \; Z_2 \mapsto \frac{1}{Z_1}, \;
Z_3 \mapsto \frac{1}{Z_3}. \\ \\
\overline{G}_{32} : & \left\{\begin{array}{ll} \displaystyle
\sigma_1 \;:\; Z_1 \mapsto Z_1, \; Z_2 \mapsto \frac{1}{Z_2}, \; Z_3
\mapsto \frac{1}{Z_3}; & \;\; \displaystyle \sigma_2 \;:\; Z_1
\mapsto \frac{1}{Z_1},
\; Z_2 \mapsto \frac{1}{Z_2}, \; Z_3 \mapsto Z_3; \\
\displaystyle \kappa \;:\; Z_1 \mapsto \frac{1}{Z_1}, \; Z_2 \mapsto
Z_2, \; Z_3 \mapsto \frac{1}{Z_3}; & \;\; \displaystyle \varphi_1
\;:\; Z_1 \mapsto
Z_3, \; Z_2 \mapsto \frac{1}{Z_1}, \; Z_3 \mapsto \frac{1}{Z_2}; \\
\psi_1^2 \;=\;
{\rm Id}. \end{array} \right. \\ \\
\overline{G}_{39} : & \; \sigma_1,\; \sigma_2,\; \kappa,\; \varphi_1
\;\; \text{as in} \;\; \overline{G}_{32}, \; \text{and} \;\;
\displaystyle \psi_1 \;:\; Z_1 \mapsto Z_1, \; Z_2 \mapsto
\frac{1}{Z_3}, \; Z_3 \mapsto \frac{1}{Z_2}. \\ \\
\overline{G}_{29} : & \left\{\begin{array}{ll} \displaystyle
\sigma_1 \;:\; Z_1 \mapsto Z_1, \; Z_2 \mapsto \frac{1}{Z_2}, \; Z_3
\mapsto
\frac{1}{Z_3}; & \; \displaystyle \sigma_2 \;=\; {\rm Id}; \\
\displaystyle \kappa \;:\; Z_1 \mapsto \frac{1}{Z_1}, \; Z_2 \mapsto
\frac{1}{Z_2}, \; Z_3 \mapsto Z_3; & \; \displaystyle \psi_2 \;:\;
Z_1 \mapsto Z_1,
\; Z_2 \mapsto \frac{1}{Z_3}, \; Z_3 \mapsto \frac{1}{Z_2}; \\
\displaystyle (2,4)(6,8) \;=\; {\rm Id}. \end{array}\right. \\ \\
\overline{G}_{22} : & \left\{\begin{array}{l} \displaystyle \sigma_1
\;=\; {\rm Id}; \quad \sigma_2 \;:\; Z_1
\mapsto Z_1, \; Z_2 \mapsto \frac{1}{Z_2}, \; Z_3 \mapsto \frac{1}{Z_3}; \\
\displaystyle \kappa \;:\; Z_1 \mapsto \frac{1}{Z_1}, \; Z_2 \mapsto
\frac{1}{Z_2}, \; Z_3 \mapsto Z_3; \quad \widetilde\varphi_2 \;=\;
{\rm Id}; \quad (3,4)(7,8) \;=\; {\rm Id}.
\end{array} \right.
\end{array}
$$
So Theorem 2.4 (2) applies.

Other 6 groups give rise to reduced monomial actions, but not pure,
so we need to check that Theorem 2.5 is applicable. In the
following, $G_{i,j,k} \subset {\rm GL}(3,\mathbb{Z})$ will be the
groups defined in \cite{HoKiYa}, section 2 (we will recall the
definitions when needed).

\bigskip

3. The groups $\overline{G}_{26}, \overline{G}_{28}$ and
$\overline{G}_{30}$ have similar behavior, in fact
$$
\begin{array}{ll}
\overline{G}_{26} : & \left\{\begin{array}{l} \displaystyle \qquad
\quad \Theta \;:\; Z_1 \mapsto -Z_2, \; Z_2 \mapsto -\frac{1}{Z_1},
\; Z_3 \mapsto -\frac{1}{Z_3}; \quad (1,5)(2,6) \;=\; {\rm Id}; \\
\displaystyle (1,5)(2,8)(4,6) \;:\; Z_1 \mapsto -\frac{1}{Z_2}, \;
Z_2 \mapsto -\frac{1}{Z_1}, \; Z_3 \mapsto -Z_3.
\end{array} \right. \\ \\
\overline{G}_{28} : & \left\{\begin{array}{l} \displaystyle
(3,7)(4,8) \;=\; {\rm Id}; \quad (2,4)(6,8) \;:\; Z_1 \mapsto
\frac{1}{Z_2}, \; Z_2 \mapsto \frac{1}{Z_1}, \; Z_3 \mapsto Z_3;
\\ \displaystyle \qquad
\quad \Theta \;\; \text{as in} \;\; \overline{G}_{26}.
\end{array} \right. \\ \\
\overline{G}_{30} : & \left\{\begin{array}{l} \displaystyle
(3,7)(4,8) \;=\; {\rm Id}; \quad (1,5)(2,4)(6,8) \;:\; Z_1 \mapsto
-\frac{1}{Z_2}, \; Z_2 \mapsto -\frac{1}{Z_1}, \; Z_3 \mapsto -Z_3; \\
\displaystyle \qquad \quad \Psi \;:\; Z_1 \mapsto Z_2, \; Z_2
\mapsto \frac{1}{Z_1}, \; Z_3 \mapsto \frac{1}{Z_3}.
\end{array} \right.
\end{array}
$$
These all correspond to $G_{4,6,1}$, so Theorem 2.5 applies. (Recall
that $G_{4,6,1}$ is defined to be $\langle -\sigma_{4A}, \lambda_1
\rangle$ with $$ \sigma_{4A} = \begin{pmatrix} 0 & -1 & 0 \\ 1 & 0 & 0 \\
0 & 0 & 1
\end{pmatrix} ,\quad \lambda_1 = \begin{pmatrix} -1 & 0 & 0 \\ 0 & 1 & 0 \\
0 & 0 & -1
\end{pmatrix}. $$
And $(Z_1 \mapsto \pm Z_2, Z_2 \mapsto \pm \frac{1}{Z_1}, Z_3
\mapsto \pm \frac{1}{Z_3})$ corresponds to $(-\sigma_{4A})^3$, while
$(Z_1 \mapsto \pm \frac{1}{Z_2}, Z_2 \mapsto \pm \frac{1}{Z_1}, Z_3
\mapsto \pm Z_3)$ corresponds to $(-\sigma_{4A})^3\lambda_1$.)

\bigskip

4. In the same manner, we check that the group $\overline{G}_{40}$
corresponds to $G_{7,4,1}$. It is better to choose new invariants
$$ Z_1' = \frac{z_3z_1}{z_2}, \quad Z_2' = \frac{z_1}{z_2z_3}, \quad
Z_3' = \frac{z_1z_2}{z_3},
$$ and we find that
$$ \overline{G}_{40} : \quad
\left\{\begin{array}{l} \displaystyle (1,5)(2,6) \;=\; {\rm Id};
\quad \sigma_1 \;:\; Z_1' \mapsto \frac{1}{Z_1'}, \;
Z_2' \mapsto \frac{1}{Z_2'}, \; Z_3' \mapsto Z_3'; \\
\displaystyle \quad \widetilde\varphi_1 \;:\; Z_1' \mapsto Z_2', \;
Z_2' \mapsto Z_3', \; Z_3' \mapsto Z_1'; \\ \displaystyle
(1,5)(3,4)(7,8) \;:\; Z_1' \mapsto -Z_2', \; Z_2' \mapsto -Z_1', \;
Z_3' \mapsto -Z_3'.
\end{array} \right.
$$ It corresponds to $G_{7,4,1}$ is now clear, since $$ G_{7,4,1} = \left\langle\begin{pmatrix} -1 & 0 & 0 \\ 0 & -1 & 0 \\ 0 & 0 & 1 \end{pmatrix} ,
\begin{pmatrix} -1 & 0 & 0 \\ 0 & 1 & 0 \\ 0 & 0 & -1 \end{pmatrix} ,
\begin{pmatrix} 0 & 0 & 1 \\ 1 & 0 & 0 \\ 0 & 1 & 0 \end{pmatrix} ,
\begin{pmatrix} 0 & -1 & 0 \\ -1 & 0 & 0 \\ 0 & 0 & 1 \end{pmatrix} \right\rangle. $$

\bigskip

5. For the groups $\overline{G}_{16}$ and $\overline{G}_{21}$,
unfortunately the corresponding groups are $G_{4,2,1}$ and
$G_{3,1,1}$, which fall into exceptional set of Theorem 2.5. One may
consult \cite{Ya} for informations, but we will give here a direct
proof of the original problems for $G_{16}$ and $G_{21}$.

Go back to $y_1,y_2,y_3,y_4$, consider the normal subgroup
$\Lambda_0 = \langle (1,5)(3,7), (2,6)(4,8)\rangle$:
$$
\begin{array}{l}
(1,5)(3,7) \;:\; y_1 \mapsto -y_1, \quad y_2 \mapsto y_2, \quad y_3 \mapsto -y_3, \quad y_4 \mapsto y_4, \\
(2,6)(4,8) \;:\; y_1 \mapsto y_1, \quad y_2 \mapsto -y_2, \quad y_3
\mapsto y_3, \quad y_4 \mapsto -y_4.
\end{array}
$$
So we have $K(y_1,y_2,y_3,y_4)^{\Lambda_0} = K(u_1,u_2,u_3,u_4)$
with $$ u_1 = y_1y_3, \quad u_2 = \frac{y_1}{y_3}, \quad u_3 =
y_2y_4, \quad u_4 = \frac{y_2}{y_4}. $$ Now let $\overline{G}_{16} =
G_{16}/\Lambda_0 \simeq C_4 \times C_2$ and $\overline{G}_{21} =
G_{21}/\Lambda_0 \simeq D_4$, we have
$$
\begin{array}{ll}
\overline{G}_{16} : & \left\{\begin{array}{l} \displaystyle \Theta
\;:\; u_1 \mapsto u_3 \mapsto -u_1,\; u_2 \mapsto u_4 \mapsto
-\frac1{u_2}; \\ \displaystyle (3,7)(4,8) \;:\; u_i \mapsto -u_i
\;\; \text{for} \;\; i=1,2,3,4;
\end{array} \right. \\ \\
\overline{G}_{21} : & \left\{\begin{array}{l} \displaystyle
(1,2,5,6)(3,4)(7,8) \;:\; u_1 \mapsto u_3 \mapsto -u_1,\; u_2
\mapsto u_4 \mapsto -u_2; \\ \displaystyle \sigma_2 \;:\; u_1
\mapsto u_1,\; u_2 \mapsto \frac1{u_2},\; u_3 \mapsto u_3,\; u_4
\mapsto \frac1{u_4}.
\end{array} \right.
\end{array}
$$
We see that $\overline{G}_{16}$ and $\overline{G}_{21}$ act on
$K(u_2,u_4)$ faithfully, so by Theorem 2.1 (1), we are reduced to
prove the rationalities for $K(u_2,u_4)^{\overline{G}_{16}}$ and
$K(u_2,u_4)^{\overline{G}_{21}}$. Now the results follows from
Theorem 2.4 (1).

\bigskip

6. There remains to consider the groups $G_{15}$ and $G_{19}$. The
group $G_{15}$ has a normal subgroup of order 4:
$$ \Lambda_5 = \langle (1,5)(3,7),\, (2,6)(4,8)\rangle , $$
and we have $K(z_1,z_2,z_3)^{\Lambda_5} = K(Z_1,Z_2,Z_3)$ with
$$ Z_1 = \frac{z_1}{z_3}, \quad Z_2 = z_2, \quad Z_3 \mapsto \frac{z_1z_3}{z_2} $$
(in fact $(2,6)(4,8) : z_1 \mapsto -z_1, z_2 \mapsto z_2, z_3
\mapsto -z_3$). The action of $\overline{G}_{15} = G_{15}/\Lambda_5$
on $K(Z_1,Z_2,Z_3)$ is monomial:
$$
\overline{G}_{15} : \quad \left\{\begin{array}{l} \displaystyle
\quad \Theta \;:\; Z_1 \mapsto Z_2, \; Z_2 \mapsto -\frac{1}{Z_1},
\; Z_3 \mapsto -\frac{1}{Z_3}; \quad (2,6)(4,8)
\;=\; {\rm Id}; \\
\displaystyle (1,8)(2,7)(3,6)(4,5) \;:\; Z_1 \mapsto \frac{1}{Z_2},
\; Z_2 \mapsto \frac{1}{Z_1}, \; Z_3 \mapsto \frac{1}{Z_3}.
\end{array} \right.
$$ This corresponds to $G_{4,6,2}$, recall that $G_{4,6,2} = \langle -\sigma_{4A},
-\lambda_1\rangle$ with $$ \sigma_{4A} = \begin{pmatrix} 0 & -1 & 0 \\ 1 & 0 & 0 \\
0 & 0 & 1
\end{pmatrix} ,\quad \lambda_1 = \begin{pmatrix} -1 & 0 & 0 \\ 0 & 1 & 0 \\
0 & 0 & -1
\end{pmatrix}, $$ and one checks that $(Z_1 \mapsto Z_2, Z_2
\mapsto -\frac{1}{Z_1}, Z_3 \mapsto -\frac{1}{Z_3})$ corresponds to
$(-\sigma_{4A})^3$, while $(Z_1 \mapsto \frac{1}{Z_2}, Z_2 \mapsto
\frac{1}{Z_1}, Z_3 \mapsto \frac{1}{Z_3})$ corresponds to
$(-\sigma_{4A})(-\lambda_1)$.

For the group $G_{19}$, it has a normal subgroup of order 4:
$$
\Lambda_6 = \langle (1,3)(2,4),\, (5,7)(6,8)\rangle ,
$$
and $K(z_1,z_2,z_3)^{\Lambda_6} = K(Z_1,Z_2,Z_3)$ with
$$ Z_1 = z_3, \quad Z_2 = \frac{z_1}{z_2}, \quad Z_3 \mapsto \frac{z_1z_2}{z_3} $$
(because $\psi_2^2 = (5,7)(6,8) : z_1 \mapsto -z_1, z_2 \mapsto
-z_2, z_3 \mapsto z_3$). The action of $\overline{G}_{19} =
G_{19}/\Lambda_6$ on $K(Z_1,Z_2,Z_3)$ is monomial:
$$
\overline{G}_{19} : \quad \left\{\begin{array}{l} \displaystyle
\sigma_1 \;:\; Z_1 \mapsto \frac{1}{Z_1}, \; Z_2 \mapsto
\frac{1}{Z_2}, \; Z_3
\mapsto Z_3; \quad \sigma_2 \;=\; {\rm Id}; \\
\displaystyle \kappa \;:\; Z_1 \mapsto Z_2, \; Z_2 \mapsto Z_1, \;
Z_3 \mapsto \frac{1}{Z_3}; \quad \psi_2 \;:\; Z_1 \mapsto
-\frac{1}{Z_1}, \; Z_2 \mapsto -Z_2, \; Z_3 \mapsto Z_3.
\end{array} \right.
$$
The corresponding group is again $G_{4,6,2}$, since $(Z_1 \mapsto
\frac{1}{Z_1}, Z_2 \mapsto \frac{1}{Z_2}, Z_3 \mapsto Z_3)$
corresponds to $(-\sigma_{4A})^2$, $(Z_1 \mapsto Z_2, Z_2 \mapsto
Z_1, Z_3 \mapsto \frac{1}{Z_3})$ corresponds to
$(-\sigma_{4A})^3(-\lambda_1)$, and $(Z_1 \mapsto -\frac{1}{Z_1},
Z_2 \mapsto -Z_2, Z_3 \mapsto Z_3)$ corresponds to
$(-\sigma_{4A})^2(-\lambda_1)$.

\bigskip

The proofs in case ${\rm char}(K) \ne 2$ is now finished.

\bigskip

Let's now consider the case ${\rm char}(K) = 2$. This is not covered
by \cite{HoKiYa}'s results, and we are even in trouble on finding
faithful subrepresentations.

Nevertheless, all these groups are 2-groups except $G_{13}, G_{14},
G_{24}, G_{32}, G_{39}, G_{40}$. Therefore in view of Theorem 2.8,
we only need to consider these 6 groups.

Among these groups, $G_{13}, G_{24}, G_{32}, G_{39}$ can be
considered together, so we divide the proof into 3 cases.

\bigskip

1. The four groups $G_{13}, G_{24}, G_{32}, G_{39}$ all have $V_8 =
\langle \sigma_1, \sigma_2, \kappa\rangle$ as a normal subgroup, and
by Proposition 2.9 (2) we have $K(x_1,\dots,x_n)^{V_8} =
K(z_1,\dots,z_8)$, where
$$ \begin{array}{l}
z_1 = x_1+x_2+x_3+x_4+x_5+x_6+x_7+x_8, \quad z_2 =
x_1x_2+x_3x_4+x_5x_6+x_7x_8, \\
z_3 = x_1x_3+x_2x_4+x_5x_7+x_6x_8, \quad
z_4 = x_1x_4+x_2x_3+x_5x_8+x_6x_7, \\
z_5 = x_1x_5+x_2x_6+x_3x_7+x_4x_8, \quad z_6 =
x_1x_6+x_2x_5+x_3x_8+x_4x_7, \\
z_7 = x_1x_7+x_2x_8+x_3x_5+x_4x_6, \quad z_8 =
x_1x_8+x_2x_7+x_3x_6+x_4x_5. \\
\end{array} $$
It is easily checked that the actions of $\varphi_1,
\widetilde\varphi_2, (3,6)(4,5)$ and $\psi_1$ on $z_i$'s are exactly
the same way as their actions on $x_i$'s, for example $\varphi_1$ on
$z_i$'s is $z_2 \mapsto z_3 \mapsto z_4 \mapsto z_2, \ z_7 \mapsto
z_6 \mapsto z_5 \mapsto z_7$ and $z_1,z_8$ fixed. This realizes
$G_{13}/V_8, G_{24}/V_8, G_{32}/V_8, G_{39}/V_8$ as subgroups of
$S_6$ (acting on $z_2,z_3,z_4,z_5,z_6,z_7$ by permutations). Now
Theorem 1.1 (2) applies (in fact, $G_{13}/V_8$ and $G_{24}/V_8$ are
nontransitive subgroups of $S_6$, isomorphic to $C_3$ and $S_3$,
while $G_{32}/V_8$ and $G_{39}/V_8$ are transitive subgroups of
$S_6$ correspond, under the notations of \cite{KWZ}, section 3, to
$G_7$ and $G_5$. Note that ${\rm PGL}(2,5)$ and ${\rm PSL}(2,5)$ are
numbered as $G_{13}$ and $G_{14}$ there).

\bigskip

2. $G_{40}$ has a normal subgroup $\Lambda_1 = A_8 \cap \langle
(1,5), (2,6), (3,7), (4,8)\rangle$, and we have
$K(x_1,\dots,x_8)^{\Lambda_1} = K(z_1,\dots,z_8)$, where $$
\begin{array}{l} z_1 = x_1+x_5, \quad z_2 = x_2+x_6, \quad z_3 =
x_3+x_7, \quad z_4 = x_4+x_8, \\ z_5 = x_1x_5, \qquad\, z_6 =
x_2x_6, \qquad\; z_7 = x_3x_7, \\ z_8 =
(x_1x_2x_3+x_1x_6x_7+x_5x_2x_7+x_5x_6x_3)x_4+
(x_5x_2x_3+x_5x_6x_7+x_1x_2x_7+x_1x_6x_3)x_8.
\end{array} $$ To see this, one checks that these are clearly
invariants of $\Lambda_1$, and the degree of extension $$
\begin{array}{rl}
& [K(x_1,\dots,x_8):K(z_1,\dots,z_8)] \\
=& [K(x_1,x_5,x_2,x_6,x_3,x_7)(z_4,z_8):K(z_1,\dots,z_8)] \\
=& [K(x_1,x_5,x_2,x_6,x_3,x_7):K(z_1,z_5,z_2,z_6,z_3,z_7)] \\
\le& [K(x_1,x_5):K(z_1,z_5)] \cdot [K(x_2,x_6):K(z_2,z_6)] \cdot
[K(x_3,x_7):K(z_3,z_7)]
\\ =& 8.
\end{array} $$

Now the actions of $\sigma_1, \widetilde\varphi_1$ and
$(1,5)(3,4)(7,8)$ on $z_i$'s are as follows:
$$ \begin{array}{rll} \sigma_1 &:& \displaystyle z_1
\leftrightarrow z_2, \quad z_3 \leftrightarrow z_4, \quad z_5
\leftrightarrow z_6, \quad z_8 \;\text{fixed}, \\ && \displaystyle
z_7 \mapsto \frac{z_4^2}{z_1^2} z_5 + \frac{z_4^2}{z_2^2} z_6 +
\frac{z_4^2}{z_3^2} z_7 + \frac{z_4}{z_1z_2z_3} z_8 +
\frac{1}{z_1^2z_2^2z_3^2} z_8^2, \\ \\
\widetilde\varphi_1 &:& \displaystyle z_2 \mapsto z_3 \mapsto z_4
\mapsto z_2, \quad z_1, z_5, z_8 \;\text{fixed}, \\ && \displaystyle
z_6 \mapsto z_7 \mapsto \frac{z_4^2}{z_1^2} z_5 +
\frac{z_4^2}{z_2^2} z_6 + \frac{z_4^2}{z_3^2} z_7 +
\frac{z_4}{z_1z_2z_3} z_8 +
\frac{1}{z_1^2z_2^2z_3^2} z_8^2, \\ \\
\hspace{-30pt}(1,5)(3,4)(7,8) &:& \displaystyle z_3 \leftrightarrow
z_4, \quad z_8 \mapsto
z_1z_2z_3z_4 + z_8, \quad z_1, z_2, z_5, z_6 \;\text{fixed}, \\
&& \displaystyle z_7 \mapsto \frac{z_4^2}{z_1^2} z_5 +
\frac{z_4^2}{z_2^2} z_6 + \frac{z_4^2}{z_3^2} z_7 +
\frac{z_4}{z_1z_2z_3} z_8 + \frac{1}{z_1^2z_2^2z_3^2} z_8^2.
\end{array} $$
Since $|G_{40}/\Lambda_1| = 24$, we find that $G_{40}/\Lambda_1$
acts on $K(z_1,z_2,z_3,z_4)$ faithfully by permutations of
$z_1,z_2,z_3,z_4$ ($\simeq S_4$), therefore
$$
\begin{array}{lll}
& K(x_1,\dots,x_8)^{G_{40}} \\
=& K(z_1,\dots,z_8)^{G_{40}/\Lambda_1} \\
=& K(z_1,z_2,z_3,z_4,z_8)^{G_{40}/\Lambda_1}(w_5,w_6,w_7) & \qquad \text{by Theorem 2.1 (1)} \\
=& K(z_1,z_2,z_3,z_4)^{G_{40}/\Lambda_1}(w_8)(w_5,w_6,w_7) & \qquad \text{by Theorem 2.1 (1) or (2)} \\
=& K(w_1,w_2,w_3,w_4)(w_8)(w_5,w_6,w_7) \\
=& K(w_1,\dots,w_8),
\end{array} $$
where $w_1,w_2,w_3,w_4$ are the elementary symmetric polynomials of
$z_1,z_2,z_3,z_4$. This finished the rationality proof for $G_{40}$.

\bigskip

3. Finally $G_{14}$ has a normal subgroup $N_4 = \langle \sigma_1,
\sigma_2\rangle \simeq V_4$. Let
$$
\begin{array}{ll} z_5 = x_5+x_6+x_7+x_8, & z_6 =
x_2x_5+x_1x_6+x_4x_7+x_3x_8, \\ z_7 = x_3x_5+x_4x_6+x_1x_7+x_2x_8, &
z_8 = x_4x_5+x_3x_6+x_2x_7+x_1x_8,
\end{array} $$ we have $K(x_1,\dots,x_8) =
K(x_1,x_2,x_3,x_4,z_5,z_6,z_7,z_8)$, and $z_5,z_6,z_7,z_8$ are fixed
by $N_4$, so by Proposition 2.9 (1), $$ K(x_1,\dots,x_8)^{N_4} =
K(x_1,x_2,x_3,x_4,z_5)^G(z_5,z_6,z_7,z_8) = K(z_1,\dots,z_8), $$
where $z_1 = x_1+x_2+x_3+x_4, \ z_2 = x_1x_2+x_3x_4, \ z_3 =
x_1x_3+x_2x_4, \ z_4 = x_1x_4+x_2x_4$. The actions of
$\widetilde\varphi_1$ and $\kappa^\circ$ on $z_i$'s are as follows:
$$ \begin{array}{lll} \widetilde\varphi_1 &:& \displaystyle z_2 \mapsto z_3
\mapsto z_4 \mapsto z_2, \quad z_6 \mapsto z_7 \mapsto z_8 \mapsto
z_6, \; \text{and}\; z_1,z_5 \;\text{fixed}, \\
\kappa^\circ  &:& \displaystyle z_1 \leftrightarrow z_5, \quad z_7
\leftrightarrow z_8, \quad z_6 \;\text{fixed}, \\
&& \displaystyle z_2 \mapsto
\frac{z_5^2}{z_1^2}z_2+\frac{z_5}{z_1}z_6+\frac{1}{z_1^2}z_6^2+A,
\quad z_3 \mapsto
\frac{z_5^2}{z_1^2}z_4+\frac{z_5}{z_1}z_8+\frac{1}{z_1^2}z_8^2+A, \\
&& \displaystyle z_4 \mapsto
\frac{z_5^2}{z_1^2}z_3+\frac{z_5}{z_1}z_7+\frac{1}{z_1^2}z_7^2+A,
\quad \text{where}\; A = \frac{z_6z_7+z_7z_8+z_8z_6}{z_1^2}.
\end{array} $$
Fix a primitive 3th root of unity $\zeta_3$ in the algebraic closure
of $K$. When $\zeta_3 \notin K$, the field $K(\zeta_3)$ is a
quadratic extension of $K$ with Galois group $\pi = \{1,\rho\}$,
where $\rho : \zeta_3 \mapsto \zeta_3^2$. In case $\zeta_3 \in K$,
we take $\pi$ to be the trivial group. Now let $$
\begin{array}{llll} w_1 = z_1, & w_2 = \zeta_3^2 z_2 + \zeta_3 z_3 +
z_4, & w_3 = \zeta_3 z_2 + \zeta_3^2 z_3 + z_4, & w_4 = z_2 + z_3 +
z_4, \\ w_5 = z_5, & w_6 = \zeta_3^2 z_6 + \zeta_3 z_7 + z_8, & w_7
= \zeta_3 z_6 + \zeta_3^2 z_7 + z_8, & w_8 = z_6 + z_7 + z_8,
\end{array} $$
the action of $\widetilde\varphi_1$ on $w_i$'s has the diagonal form
${\rm diag}(1,\zeta_3,\zeta_3^2,1,1,\zeta_3,\zeta_3^2,1)$, and that
of $\kappa^\circ$ becomes:
$$ \begin{array}{l} \displaystyle w_1 \leftrightarrow w_5, \quad w_6 \mapsto
\zeta_3 w_7, \quad w_7 \mapsto \zeta_3^2 w_6, \quad w_8
\;\text{fixed}, \\ \displaystyle w_2 \mapsto \zeta_3
\Big(\frac{w_5^2}{w_1^2} w_3 + \frac{w_5}{w_1} w_7 + \frac{1}{w_1^2}
w_6^2\Big), \quad w_3 \mapsto \zeta_3^2 \Big(\frac{w_5^2}{w_1^2} w_2
+ \frac{w_5}{w_1} w_6 + \frac{1}{w_1^2} w_7^2\Big), \\ \displaystyle
w_4 \mapsto \frac{w_5^2}{w_1^2} w_4 + \frac{w_5}{w_1} w_8 +
\frac{1}{w_1^2} w_6w_7.
\end{array} $$
Note that if $\rho$ exists, it acts as: $$ w_2 \leftrightarrow w_3,
\quad w_6 \leftrightarrow w_7, \; \text{and} \; w_1,w_4,w_5,w_8
\;\text{fixed}. $$ Let $$ u_2 = \frac{w_5}{w_1w_7^2}w_2, \quad u_3 =
\frac{w_5}{w_1w_6^2}w_3, \quad u_4 = \frac{w_5}{w_1}w_4,
$$ we have $K(\zeta_3)(w_1,\dots,w_8) =
K(\zeta_3)(w_1,w_5,w_6,w_7,w_8)(u_2,u_3,u_4)$, and
$u_2,u_3,u_4,w_1,w_5,$ $w_8$ are all fixed by $\widetilde\varphi_1$.
It follows that the fixed field of $\widetilde\varphi_1$ has the
form $K(\zeta_3)(u_1,\dots,u_8)$, where
$$ u_1 = w_1, \quad u_5 = w_5, \quad u_6 = \frac{w_6}{w_7^2}, \quad u_7 = \frac{w_7}{w_6^2}, \quad u_8
= w_8.
$$ The choice of $u_i$'s makes the action of $\kappa^\circ$
into the following simple form: $$ \begin{array}{l} \displaystyle
u_1 \leftrightarrow u_5, \quad u_6 \leftrightarrow u_7, \quad u_8
\;\text{fixed}, \\ \displaystyle u_2 \mapsto u_3 + u_7 +
\frac1{u_1u_5}, \quad u_3 \mapsto u_2 + u_6 + \frac1{u_1u_5}, \quad
u_4 \mapsto u_4 + u_8 + \frac1{u_1u_5u_6u_7}. \end{array}
$$
Finally, let $$ \begin{array}{l} \displaystyle v_1 = \zeta_3^2 u_1 +
\zeta_3 u_5, \quad v_5 = \zeta_3 u_1 + \zeta_3^2 u_5, \quad
v_6 = u_6, \quad v_7 = u_7, \quad v_8 = u_8, \\
\displaystyle v_2 = u_2 +
\frac{u_1}{u_1+u_5}\Big(u_6+\frac1{u_1u_5}\Big), \quad v_3 = u_3 +
\frac{u_1}{u_1+u_5}\Big(u_7+\frac1{u_1u_5}\Big), \\ \displaystyle
v_4 = u_4 + \frac{u_1}{u_1+u_5}\Big(u_8+\frac1{u_1u_5u_6u_7}\Big).
\end{array} $$ One sees that \begin{multline*}
K(\zeta_3)(u_1,\dots,u_8) =
K(\zeta_3)(v_1,v_5,v_6,v_7,v_8)(u_2,u_3,u_4) \\ =
K(\zeta_3)(v_1,v_5,v_6,v_7,v_8)(v_2,v_3,v_4) =
K(\zeta_3)(v_1,\dots,v_8),
\end{multline*} and computations show that $\kappa^\circ$
acts on $v_i$'s as follows:
$$ v_1 \leftrightarrow v_5, \quad v_2 \leftrightarrow v_3, \quad v_6
\leftrightarrow v_7, \quad \text{and}\; v_4, v_8 \;\text{fixed}. $$
If $\zeta_3 \notin K$, then the abelian group $\langle
\kappa^\circ\rangle \times \pi$ is generated by $\kappa^\circ\rho$
and $\kappa^\circ$, and one checks immediately that $v_1,\dots,v_8$
are fixed by $\kappa^\circ\rho$, and $\zeta_3$ mapped to
$\zeta_3^2$, this shows $K(\zeta_3)(v_1,\dots,v_8)^{\langle
\kappa^\circ\rho\rangle} = K(\zeta_3)^{\langle
\kappa^\circ\rho\rangle}(v_1,\dots,v_8) = K(v_1,\dots,v_8)$. So in
any case, the subfield of $K(\zeta_3)(v_1,\dots,v_8)$ fixed by
$\langle \kappa^\circ\rangle \times \pi$ is
$$ K(v_1,\dots,v_8)^{\langle \kappa^\circ\rangle} = K(v_1+v_5,
v_2+v_3,v_1v_3+v_5v_2,v_4,v_1v_5,v_6+v_7,v_1v_7+v_5v_6,v_8) $$ by
Proposition 2.2, so it is $K$-rational.

\bigskip

For $G_{14}$ and $G_{40}$, there is another approach to prove the
rationality in case of characteristic 2, see the Remark at the end
of \S4.

\section{Subgroups related to $S_6$ actions.}

These are the groups $G_{33}$, $G_{34}$, $G_{41}$, $G_{45}$ and
$G_{46}$, the following group $$ V_4 \times V_4 = \langle
(1,2)(3,4),\; (1,3)(2,4),\; (5,6)(7,8),\; (5,7)(6,8)\rangle $$ is a
normal subgroup for all these groups.

First assume ${\rm char}(K) \ne 2$, by taking variable change
$$
\begin{array}{ll}
y_1 = x_1 + x_2 + x_3 + x_4, & y_2 = x_1 + x_2 - x_3 - x_4, \\
y_3 = x_1 - x_2 + x_3 - x_4, & y_4 = x_1 - x_2 - x_3 + x_4, \\
y_5 = x_5 + x_6 + x_7 + x_8, & y_6 = x_5 + x_6 - x_7 - x_8, \\
y_7 = x_5 - x_6 + x_7 - x_8, & y_8 = x_5 - x_6 - x_7 + x_8,
\end{array} $$
we find that the action of $V_4 \times V_4$ on $y_i$'s becomes:
$$ \begin{array}{lll}
(1,2)(3,4) &=& {\rm diag}(1,1,-1,-1,1,1,1,1), \\
(1,3)(2,4) &=& {\rm diag}(1,-1,1,-1,1,1,1,1), \\
(5,6)(7,8) &=& {\rm diag}(1,1,1,1,1,1,-1,-1), \\
(5,7)(6,8) &=& {\rm diag}(1,1,1,1,1,-1,1,-1).
\end{array}
$$
Using Lemma 2.11, one checks that $K(y_1,\dots,y_8)^{V_4 \times V_4}
= K(z_1,\dots,z_8)$, with
$$ z_1 = \frac{y_2 y_3}{y_4},\; z_2
= \frac{y_3 y_4}{y_2},\; z_3 = \frac{y_4 y_2}{y_3},\; z_4 =
\frac{y_6 y_7}{y_8},\; z_5 = \frac{y_7 y_8}{y_6},\; z_6 = \frac{y_8
y_6}{y_7},\; z_7 = y_1,\; z_8 = y_5.
$$

The subfield $K(z_1,\dots,z_6)$ is stable under these $G_i$'s, and
the actions of various elements on $z_1,\dots,z_6$'s are as follows:
$$
\begin{array}{rll}
\kappa &:& z_1 \leftrightarrow z_4,\; z_2 \leftrightarrow z_5,\; z_3
\leftrightarrow z_6,
\\ \widetilde\kappa  &:& z_1 \leftrightarrow z_6,\; z_2 \leftrightarrow z_5,\;
z_3 \leftrightarrow z_4,
\\ \kappa' &:& z_1 \mapsto z_6 \mapsto z_3 \mapsto z_4 \mapsto z_1,\;
z_2 \leftrightarrow z_5,
\\ \varphi_1 &:& z_1 \mapsto z_2 \mapsto z_3 \mapsto z_1,\; z_4 \mapsto
z_5 \mapsto z_6 \mapsto z_4,
\\ \widetilde\varphi_1 &:& z_1 \mapsto z_2 \mapsto z_3 \mapsto z_1,\; z_4
\mapsto z_5 \mapsto z_6 \mapsto z_4,
\\ \psi_2 &:& z_1 \leftrightarrow z_2,\; z_4 \leftrightarrow z_5,\;
\text{and}\; z_3, z_6 \;\text{are fixed},
\\ (2,3,4) &:& z_1 \mapsto z_2 \mapsto z_3 \mapsto z_1,\; \text{and}\;
z_4, z_5, z_6 \;\text{are fixed},
\\
(1,2)(5,6) &:& z_1 \leftrightarrow z_3,\; z_4 \leftrightarrow z_6,\;
\text{and}\; z_2, z_5 \;\text{are fixed}.
\end{array}
$$
So the actions of these $G_i$'s on $z_1,\dots,z_6$ are all by
permutations. By writing $\overline{G}_i \;=\; G_i/(V_4 \times
V_4)$, one verifies without difficulty (by comparing generators and
counting orders) that the actions of
$\overline{G}_{33},\overline{G}_{34},\overline{G}_{41},\overline{G}_{45},
\overline{G}_{46}$ on $z_1,\dots,z_6$ realize them as transitive
subgroups of $S_6$ correspond respectively to $G_1, G_2, G_3,
G_{11}, G_{10}$ under the notations of \cite{KWZ}, section 3
(remember that ${\rm PGL}(2,5)$ and ${\rm PSL}(2,5)$ are numbered as
$G_{13}$ and $G_{14}$ there). This shows also that the actions of
these $G_i$ on $z_1,\dots,z_6$ are all faithful, so by Theorem 2.1
(1),
$$ K(x_1,\dots,x_8)^{G_i} = K(z_1,\dots,z_6,z_7,z_8)^{\over{G}_i} =
K(z_1,\dots,z_6)^{\over{G}_i}(w_7,w_8). $$ Then by Theorem 1.1 (2),
$K(z_1,\dots,z_6)^{\over{G}_i} = K(w_1,\dots,w_6)$, and we are done.

\bigskip

Now we consider the cases when ${\rm char}(K)=2$. By Proposition 2.9
(1), the fixed field $K(x_1,\dots,x_8)^{V_4 \times V_4}$ is
generated by
$$
\begin{array}{l} z_1 = x_1+x_2+x_3+x_4, \\
z_2 = x_1x_2+x_3x_4, \quad z_3 = x_1x_3+x_2x_4, \quad z_4 =
x_1x_4+x_2x_3, \\
z_5 = x_5+x_6+x_7+x_8, \\
z_6 = x_5x_6+x_7x_8, \quad z_7 = x_5x_7+x_6x_8, \quad z_8 =
x_5x_8+x_6x_7.
\end{array} $$
We divide the 5 groups into 3 sets, and consider them separately.

\bigskip

1. For $G_{45}$ and $G_{46}$, a more bigger normal subgroup exists,
which is $A_4 \times A_4$, with the first $A_4$ (resp. the second
$A_4$) acting by permutations on $x_1,x_2,x_3,x_4$ (resp.
$x_5,x_6,x_7,x_8$). As $A_4/V_4 \simeq C_3$, we have by Theorem 2.10
$K(x_1,\dots,x_8)^{A_4 \times A_4} = K(z_1,\dots,z_8)^{C_3 \times
C_3} = K(w_1,\dots,w_8)$, with $$ \begin{array}{l} \displaystyle w_1
= z_1, \quad
w_2 = z_2+z_3+z_4, \\
\displaystyle w_3 =
\frac{z_2z_3^2+z_3z_4^2+z_4z_2^2+z_2z_3z_4}{z_2^2+z_3^2+z_4^2+z_2z_3+z_3z_4+z_4z_2},
\quad w_4 =
\frac{z_2z_4^2+z_3z_2^2+z_4z_3^2+z_2z_3z_4}{z_2^2+z_3^2+z_4^2+z_2z_3+z_3z_4+z_4z_2},
\\ \displaystyle w_5 = z_5, \quad
w_6 = z_6+z_7+z_8, \\
\displaystyle w_7 =
\frac{z_6z_7^2+z_7z_8^2+z_8z_6^2+z_6z_7z_8}{z_6^2+z_7^2+z_8^2+z_6z_7+z_7z_8+z_8z_6},
\quad w_8 =
\frac{z_6z_8^2+z_7z_6^2+z_8z_7^2+z_6z_7z_8}{z_6^2+z_7^2+z_8^2+z_6z_7+z_7z_8+z_8z_6},
\end{array} $$ since the two $C_3$ acts on $z_i$'s by
$z_1 \mapsto z_1, z_2 \mapsto z_3 \mapsto z_4 \mapsto z_2$ and $z_5
\mapsto z_5, z_6 \mapsto z_7 \mapsto z_8 \mapsto z_6$ respectively.
Now $G_{45}/(A_4 \times A_4)$ and $G_{46}/(A_4 \times A_4)$ are all
isomorphic to $C_2 \times C_2$, with common first factor $C_2$
generated by $(1,2)(5,6)$. The action of $(1,2)(5,6)$ on $w_i$'s is
as follows:
$$ w_3 \leftrightarrow w_4, \quad w_7 \leftrightarrow w_8, \;
\text{and}\; w_1,w_2,w_5,w_6 \;\text{fixed}. $$ By Proposition 2.2,
we may write the fixed field of $(1,2)(5,6)$ as $K(u_1,\dots,u_8)$,
with $$
\begin{array}{l} u_1 = w_1, \quad u_2 = w_2, \quad u_3 = w_3+w_4,
\quad u_4 = w_3w_4, \\ u_5 = w_5, \quad u_6 = w_6, \quad u_7 =
w_7+w_8, \quad u_8 = w_3w_8+w_4w_7.
\end{array} $$ In case $G_{45}$ (resp. $G_{46}$), the second factor
$C_2$ is generated by $\kappa$ (resp. $\kappa'$), we check that the
action of $\kappa$ on $u_i$'s is
$$ u_1 \leftrightarrow u_5, \quad u_2 \leftrightarrow u_6, \quad u_3
\leftrightarrow u_7, \quad u_4 \mapsto \frac{u_7^2}{u_3^2}u_4 +
\frac{u_7}{u_3}u_8 + \frac{1}{u_3^2}u_8^2, \quad \text{and}\; u_8
\;\text{fixed},
$$ while the action of $\kappa'$ is almost the
same, with only one exception $u_8 \mapsto u_3u_7+u_8$. Let $u_4' =
\frac{u_7}{u_3}u_4$, we see that $\kappa$ and $\kappa'$ maps $u_4'$
to $u_4' + \frac{u_8(u_8+u_3u_7)}{u_3u_7}$. Put $v_4 =
u_4'+\frac{u_1}{u_1+u_5}\frac{u_8(u_8+u_3u_7)}{u_3u_7}$, and let
$v_8 = u_8$ for $G_{45}$ and $=u_8+\frac{u_1}{u_1+u_5}u_3u_7$ for
$G_{46}$. It is trivially verified that in each case
$K(u_1,\dots,u_8) = K(u_1,u_2,u_3,u_5,u_6,u_7)(v_4,v_8)$, and
$v_4,v_8$ are invariant under the corresponding group. Now the
result of Proposition 2.2 applies.

\bigskip

2. For $G_{34}$, let $\zeta_3$ be a primitive 3th root of unity in
the algebraic closure of $K$, and $\pi = {\rm Gal}(K(\zeta_3)/K) =
\{1\}$ or $\{1,\rho\}$. Set $$
\begin{array}{llll} w_1 = z_1, & w_2 = \zeta_3^2 z_2 + \zeta_3 z_3 +
z_4, & w_3 = \zeta_3 z_2 + \zeta_3^2 z_3 + z_4, & w_4 = z_2 + z_3 +
z_4, \\ w_5 = z_5, & w_6 = \zeta_3^2 z_6 + \zeta_3 z_7 + z_8, & w_7
= \zeta_3 z_6 + \zeta_3^2 z_7 + z_8, & w_8 = z_6 + z_7 + z_8,
\end{array} $$ we see that $$ \begin{array}{lll}
\widetilde\varphi_1 &:& w_2 \mapsto \zeta_3 w_2, \quad w_3 \mapsto
\zeta_3^2 w_3, \quad w_6 \mapsto \zeta_3 w_6, \quad w_7 \mapsto
\zeta_3^2 w_7, \\ &&
\text{and}\; w_1,w_4,w_5,w_8 \;\text{fixed}, \\
\widetilde\kappa &:& w_2 \mapsto \zeta_3 w_7, \quad w_3 \mapsto
\zeta_3^2 w_6, \quad w_6 \mapsto \zeta_3 w_3, \quad w_7 \mapsto
\zeta_3^2 w_2, \\ && w_1 \leftrightarrow w_5, \quad w_4
\leftrightarrow w_8.
\end{array} $$ Also, if $\rho$ exists, it acts on $w_i$'s as
follows:
$$ w_2 \leftrightarrow w_3, \quad w_6 \leftrightarrow w_7, \; \text{and}\;
w_1,w_4,w_5,w_8 \;\text{fixed}. $$

We may write the fixed field of $\widetilde\varphi_1$ as
$K(\zeta_3)(u_1,\dots,u_8)$, with $$ u_1 = w_1, \quad u_2 =
\frac{w_2}{w_6}, \quad u_3 = \frac{1}{w_2w_7}, \quad u_4 = w_4, $$
$$ u_5 = w_5, \quad u_6 = \frac{w_6}{w_3w_7}, \quad u_7 =
\frac{w_7}{w_3}, \quad u_8 = w_8
$$
(use Lemma 2.11 to compute the extension degree), and the action of
$\widetilde\kappa$ becomes
$$ u_1 \leftrightarrow u_5, \quad u_2 \leftrightarrow u_7, \quad
u_3 \;\text{fixed}, \quad u_4 \leftrightarrow u_8, \quad u_6 \mapsto
\frac{u_3}{u_6}. $$ Now let $$ \begin{array}{l} \displaystyle v_1 =
u_1+u_5, \quad v_5 = u_2u_5+u_7u_1, \quad v_3 = u_6+\frac{u_3}{u_6},
\\ \displaystyle v_6 = u_2u_6+u_7\frac{u_3}{u_6}, \quad v_4 = u_4+u_8, \quad v_8 =
u_2u_8+u_7u_4, \end{array}
$$ we see easily that $K(\zeta_3)(u_1,\dots,u_8) =
K(\zeta_3)(u_2,u_7)(v_1,v_3,v_4,v_5,v_6,v_8)$, and $v_1,v_3,v_4,$
$v_5,v_6,v_8$ are all fixed by $\widetilde\kappa$, so the fixed
field of $\widetilde\kappa$ is $K(\zeta_3)(v_1,\dots,v_8)$ with $v_2
= u_2+u_7, v_7 = u_2u_7$.

We have proved $$ K(\zeta_3)(x_1,\dots,x_8)^{G_{34}} =
K(\zeta_3)(v_1,\dots,v_8),
$$ for $\zeta_3 \in K$ the proof is finished, and if $\zeta_3 \notin K$,
it remains to consider the action of $\rho$. By calculation, we have
$$ \begin{array}{lll} \rho &:& \displaystyle v_3 \leftrightarrow v_6, \quad
v_1,v_4 \;\text{fixed}, \\ && \displaystyle v_7 \mapsto
\frac{1}{v_7}, \quad v_2 \mapsto \frac{v_2}{v_7}, \quad v_5 \mapsto
\frac{v_5}{v_7}, \quad v_8 \mapsto \frac{v_8}{v_7}.
\end{array}
$$
Put $$ V_1 = v_1, \quad V_4 = v_4, \quad V_3 = \zeta_3^2 v_3+\zeta_3
v_6, \quad V_6 = \zeta_3 v_3+\zeta_3^2 v_6 $$ $$ V_2 =
v_2\Big(1+\frac1{v_7}\Big), \quad V_5 = v_5\Big(1+\frac1{v_7}\Big),
\quad V_8 = v_8\Big(1+\frac1{v_7}\Big),
$$ one finds
$$ K(\zeta_3)(v_7)(v_1,v_2,v_3,v_4,v_5,v_6,v_8) = K(\zeta_3)(v_7)(V_1,V_2,V_3,V_4,V_5,V_6,V_8). $$
Let $V_7 = \frac{1}{v_7+1}+\zeta_3$, then $K(\zeta_3)(v_7) =
K(\zeta_3)(V_7)$, and all these $V_i$'s are fixed by $\rho$, so
finally $$ K(\zeta_3)(v_1,\dots,v_8)^\pi =
K(\zeta_3)(V_1,\dots,V_8)^\pi = K(V_1,\dots,V_8), $$ we are done.

We can also use the method indicated at the end of section \S4 to
prove the result for $G_{34}$ in characteristic 2 case.

\bigskip

3. For $G_{33}$ and $G_{41}$, we use again
$$
\begin{array}{llll} w_1 = z_1, & w_2 = \zeta_3^2 z_2 + \zeta_3 z_3 +
z_4, & w_3 = \zeta_3 z_2 + \zeta_3^2 z_3 + z_4, & w_4 = z_2 + z_3 +
z_4, \\ w_5 = z_5, & w_6 = \zeta_3^2 z_6 + \zeta_3 z_7 + z_8, & w_7
= \zeta_3 z_6 + \zeta_3^2 z_7 + z_8, & w_8 = z_6 + z_7 + z_8.
\end{array} $$ The actions of $\varphi_1, \kappa$
and $\psi_2$ are now
$$ \begin{array}{lll}
\varphi_1 &:& w_2 \mapsto \zeta_3 w_2, \quad w_3 \mapsto \zeta_3^2
w_3, \quad w_6 \mapsto \zeta_3 w_6, \quad w_7 \mapsto \zeta_3^2 w_7,
\\ &&
\text{and}\; w_1,w_4,w_5,w_8 \;\text{fixed}, \\
\kappa &:& w_1 \leftrightarrow w_5, \quad w_2 \leftrightarrow w_6,
\quad w_3 \leftrightarrow w_7, \quad w_4
\leftrightarrow w_8, \\
\psi_2 &:& w_2 \mapsto \zeta_3^2 w_3, \quad w_3 \mapsto \zeta_3 w_2,
\quad w_6 \mapsto \zeta_3^2 w_7, \quad w_7 \mapsto \zeta_3 w_6, \\
&& \text{and}\; w_1,w_4,w_5,w_8 \;\text{fixed}.
\end{array} $$
Use Lemma 2.11, the fixed field of $\varphi_1$ can be written as
$K(\zeta)(u_1,\dots,u_8)$, where $$ u_1 = w_1, \quad u_2 =
\frac{w_2}{w_3w_7}, \quad u_3 = \frac{w_3}{w_2w_6}, \quad u_4 = w_4,
$$ $$ u_5 = w_5, \quad u_6 = \frac{w_6}{w_3w_7}, \quad u_7 =
\frac{w_7}{w_2w_6}, \quad u_8 = w_8.
$$ One verifies that $$ \begin{array}{lll}
\kappa &:& u_1 \leftrightarrow u_5, \quad u_2 \leftrightarrow u_6,
\quad u_3 \leftrightarrow u_7, \quad u_4
\leftrightarrow u_8, \\
\psi_2 &:& u_2 \leftrightarrow u_3, \quad u_6 \leftrightarrow u_7,
\quad \text{and}\; u_1,u_4,u_5,u_8 \;\text{fixed}.
\end{array} $$
Therefore by Proposition 2.2, $$ K(\zeta_3)(x_1,\dots,x_8)^{G_{33}}
= K(\zeta_3)(u_1,\dots,u_8)^{\langle \kappa\rangle} =
K(\zeta_3)(v_1,\dots,v_8),
$$
for $v_1=u_1+u_5, v_2=u_2+u_6, v_3=u_3+u_7, v_4=u_4+u_8, v_5=u_1u_5,
v_6=u_1u_6+u_5u_2, v_7=u_1u_7+u_5u_3, v_8=u_1u_8+u_5u_4$. It is
easily checked that the action of $\psi_2$ on $v_i$'s is the same
way as on $u_i$'s, therefore $$ K(\zeta_3)(x_1,\dots,x_8)^{G_{41}} =
K(\zeta_3)(v_1,\dots,v_8)^{\langle \psi_2\rangle} =
K(\zeta_3)(t_1,\dots,t_8), $$ with $t_1=v_1, t_2=v_2+v_3,
t_3=v_2v_3, t_4=v_4, t_5=v_5, t_6=v_6+v_7, t_7=v_2v_7+v_3v_6,
t_8=v_8$.

Finally, direct computations show that $\rho$ (if exists) acts on
$v_i$'s as follows:
$$ v_2 \leftrightarrow v_3, \quad v_6
\leftrightarrow v_7, \quad \text{and}\; v_1,v_4,v_5,v_8
\;\text{fixed}, $$ and $t_i$'s are all fixed by $\rho$. So we have
$$ K(\zeta_3)(v_1,\dots,v_8)^{\langle\rho\rangle} = K(v_1,\zeta_3^2 v_2+\zeta_3 v_3,
\zeta_3 v_2+\zeta_3^2 v_3, v_4, v_5, \zeta_3^2 v_6+\zeta_3 v_7,
\zeta_3 v_6+\zeta_3^2 v_7, v_8) $$ and
$$ K(\zeta_3)(t_1,\dots,t_8)^{\langle\rho\rangle} = K(t_1,\dots,t_8), $$
the proof is completed.

\section{Subgroups related to $S_7$ actions.}

These are the groups $G_{25}$, $G_{36}$ and $G_{48}$.

All three groups contain the group $V_8 = \langle
\sigma_1,\sigma_2,\kappa\rangle \simeq C_2 \times C_2 \times C_2$ as
a normal subgroup. An element of order 7 is $\theta$.

We have $$ G_{25} \simeq V_8 \rtimes {\it\Gamma}_0 ,\qquad G_{36}
\simeq V_8 \rtimes {\it\Gamma}_1 ,\qquad G_{48} \simeq V_8 \rtimes
{\it\Gamma}_2, $$ where ${\it\Gamma}_0 = \langle \theta \rangle$,
${\it\Gamma}_1 = \langle \theta, \varphi_1\rangle$, ${\it\Gamma}_2 =
\langle \theta, \varphi_2 \rangle$. One can verify that
${\it\Gamma}_1 \simeq C_7 \rtimes C_3$, ${\it\Gamma}_2 \simeq {\rm
PSL}(2,7)$ and ${\it\Gamma}_1 < {\it\Gamma}_2$.

\bigskip

As usual, let's start with the assumption ${\rm char}(K) \ne 2$. We
will prove the rationality of $G_{25}, G_{36}$ and $G_{48}$ in three
steps. The first step is to determine the field of invariants for
$V_8$, which is done by a suitable change of variables (see below).
In the second step, we consider further those invariants for
${\it\Gamma}_0, {\it\Gamma}_1$ and ${\it\Gamma}_2$ respectively.
Here the problem is about some 7-dimensional monomial actions, and
these can be linearized to problems about subgroups of $S_7$, so
that Theorem 1.1 (3) applies. In the final step, we invoke a
homogeneity consideration to establish our results.

Step 1. Choose the following change of variables
$$
\begin{array}{lll}
y_1 &=& (x_1 + x_2 + x_3 + x_4) + (x_5 + x_6 + x_7 + x_8), \\
y_2 &=& (x_1 - x_2 - x_3 + x_4) + (x_5 - x_6 - x_7 + x_8), \\
y_3 &=& (x_1 + x_2 - x_3 - x_4) + (x_5 + x_6 - x_7 - x_8), \\
y_4 &=& (x_1 + x_2 - x_3 - x_4) - (x_5 + x_6 - x_7 - x_8), \\
y_5 &=& (x_1 + x_2 + x_3 + x_4) - (x_5 + x_6 + x_7 + x_8), \\
y_6 &=& (x_1 - x_2 + x_3 - x_4) + (x_5 - x_6 + x_7 - x_8), \\
y_7 &=& (x_1 - x_2 + x_3 - x_4) - (x_5 - x_6 + x_7 - x_8), \\
y_8 &=& (x_1 - x_2 - x_3 + x_4) - (x_5 - x_6 - x_7 + x_8).
\end{array} $$
This is the transformation that decomposes the representation of
$V_8$ into irreducible ones, the action of $V_8$ on $y_i$'s now has
the form:
$$ \begin{array}{lll} \sigma_1 &=& {\rm diag}(1,-1,1,1,1,-1,-1,-1), \\
\sigma_2 &=& {\rm diag}(1,-1,-1,-1,1,1,1,-1), \\ \kappa &=& {\rm
diag}(1,1,1,-1,-1,1,-1,-1). \end{array} $$ Using Lemma 2.11, one
checks directly that $K(x_1,\dots,x_8)^{V_8} = K(z_1,\dots,z_8)$,
where
$$ z_1 = y_1,\; z_2 = \frac{y_2
y_3}{y_6},\; z_3 = \frac{y_3 y_7}{y_8},\; z_4 = \frac{y_4
y_5}{y_3},\; z_5 = \frac{y_5 y_6}{y_7},\; z_6 = \frac{y_6
y_8}{y_4},\; z_7 = \frac{y_7 y_4}{y_2},\; z_8 = \frac{y_8 y_2}{y_5}.
$$

Step 2. The action of $\theta$ on $z_i$'s has the same form as on
$x_i$'s:
$$ \theta \;:\; z_1 \mapsto z_1,\; z_2 \mapsto z_3 \mapsto z_7 \mapsto z_4
\mapsto z_5 \mapsto z_6 \mapsto z_8 \mapsto z_2. $$ And by direct
calculations (with the help of a computer), we find the actions of
$\varphi_1$ and $\varphi_2$ are as follows:
\begin{multline*}
\varphi_1 \;:\; z_1 \mapsto z_1,\; z_2 \mapsto \frac{z_5 z_8
z_3}{z_6 z_2},\; z_3 \mapsto \frac{z_8 z_3 z_4}{z_2 z_7},\; z_4
\mapsto \frac{z_4 z_6 z_2}{z_5 z_8}, \\ z_5 \mapsto \frac{z_6 z_2
z_7}{z_8 z_3},\; z_6 \mapsto \frac{z_2 z_7 z_5}{z_3 z_4},\; z_7
\mapsto \frac{z_3 z_4 z_6}{z_7 z_5},\; z_8 \mapsto \frac{z_7 z_5
z_8}{z_4 z_6};
\end{multline*}
\begin{multline*}
\varphi_2 \;:\; z_1 \mapsto z_1,\; z_2 \mapsto \frac{z_8 z_3
z_4}{z_2 z_7},\; z_3 \mapsto \frac{z_5 z_8 z_3}{z_6 z_2},\; z_4
\mapsto \frac{z_4 z_6 z_2}{z_5 z_8}, \\ z_5 \mapsto \frac{z_6
z_2}{z_8},\; z_6 \mapsto \frac{z_2 z_7 z_5}{z_3 z_4},\; z_7 \mapsto
z_7,\; z_8 \mapsto z_8.
\end{multline*}
Therefore, $K(z_1,\dots,z_8)$ carries a monomial action of
${\it\Gamma}_2 \ (> {\it\Gamma}_1 > {\it\Gamma}_0)$.

A surprising fact is that this action is linearizable, namely we let
$$
\begin{array}{l} w_0 = z_2 z_3 z_4 z_5 z_6 z_7 z_8, \quad w_1 = z_1,
\\
w_2 = z_2 z_7 z_5, \quad w_3 = z_3 z_4 z_6, \quad w_4 = z_4 z_6 z_2, \\
w_5 = z_5 z_8 z_3, \quad w_6 = z_6 z_2 z_7, \quad w_7 = z_7 z_5 z_8,
\quad w_8 = z_8 z_3 z_4.
\end{array}
$$ One see easily that $K(z_1,\dots,z_8) = K(w_0,w_1,\dots,w_8)$,
and
$$
\begin{array}{ll} \; \theta \;\;:\; & w_2
\mapsto w_3 \mapsto w_7 \mapsto w_4 \mapsto w_5 \mapsto w_6 \mapsto
w_8 \mapsto w_2,\; \text{and}\; w_0, w_1 \;\text{fixed}; \\
\varphi_1 \;: & w_2 \mapsto w_3 \mapsto w_4 \mapsto w_2,\; w_7
\mapsto w_6 \mapsto w_5 \mapsto w_7,\; \text{and}\; w_0, w_1, w_8
\;\text{fixed}; \\ \varphi_2 \;: & w_2 \leftrightarrow w_3,\; w_7
\leftrightarrow w_6,\; \text{and}\; w_0, w_1, w_4, w_5, w_8
\;\text{fixed}.
\end{array}
$$ Note that the actions of $\theta,\varphi_1,\varphi_2$ on $w_1,\dots,w_8$
are exactly the same way as their actions on the $x_i$'s, this
explains our previous choices.

There is also a quite simple relation between these $w_i$'s, which
is
$$ w_0^3 = w_2 w_3 w_4 w_5 w_6 w_7 w_8. $$

Step 3. The above actions of ${\it\Gamma}_0, {\it\Gamma}_1$ and
${\it\Gamma}_2$ on $w_2,\dots,w_8$ make them into transitive
subgroups of $S_7$. So by Theorem 1.1 (3),
$K(w_2,\dots,w_8)^{{\it\Gamma}_0}$ and
$K(w_2,\dots,w_8)^{{\it\Gamma}_1}$ are $K$-rational, and so is
$K(w_2,\dots,w_8)^{{\it\Gamma}_2}$ provides $\mathbb{Q}(\sqrt{-7})
\subseteq K$.

But this is not enough, and to go further, we need to take $w_1$
into account. For brevity, let's put ${\it\Gamma} = {\it\Gamma}_0,
{\it\Gamma}_1$ or ${\it\Gamma}_2$, and assume $\mathbb{Q}(\sqrt{-7})
\subseteq K$ when ${\it\Gamma} = {\it\Gamma}_2$.

In view of Proposition 2.6, we have $$ K(w_1, \dots,
w_8)^{\it\Gamma} = K(v_1, \dots, v_8),
$$ where the $v_i$'s are homogeneous rational functions of $w_1,\dots,w_8$ of degree 1.
Now $w_2 w_3 w_4 w_5 w_6 w_7 w_8$ is a homogeneous element of degree
7, and belongs to $K(w_1, \dots, w_8)^{\it\Gamma}$, so it can be
written as a rational function $P(v_1, \dots, v_8)$ homogeneous of
degree 7 (see the discussion before Proposition 2.6).

Put everything together, we have shown that $$
K(z_1,\dots,z_8)^{\it\Gamma} = K(w_0,w_1,\dots,w_8)^{\it\Gamma} =
K(w_0, v_1, \dots, v_8)
$$
with one relation $w_0^3 = P(v_1, \dots, v_8)$, and $P$ homogeneous
of degree 7. Now we can complete the rationality proof by using
Lemma 2.7.

\bigskip

Finally, we come to the characteristic 2 case. By Proposition 2.9
(2), we have $K(x_1,\dots,x_8)^{V_8} = K(w_1,\dots,w_8)$ with $$
\begin{array}{l} w_1 = x_1+x_2+x_3+x_4+x_5+x_6+x_7+x_8, \quad
w_2 = x_1x_2+x_3x_4+x_5x_6+x_7x_8, \\
w_3 = x_1x_3+x_2x_4+x_5x_7+x_6x_8, \quad
w_4 = x_1x_4+x_2x_3+x_5x_8+x_6x_7, \\
w_5 = x_1x_5+x_2x_6+x_3x_7+x_4x_8, \quad
w_6 = x_1x_6+x_2x_5+x_3x_8+x_4x_7, \\
w_7 = x_1x_7+x_2x_8+x_3x_5+x_4x_6, \quad w_8 =
x_1x_8+x_2x_7+x_3x_6+x_4x_5,
\end{array} $$ and the actions of $\theta, \varphi_1, \varphi_2$
acts on $w_i$'s exactly the same way as their actions on $x_i$'s.
But then the actions of $G_{25}/V_8, G_{36}/V_8, G_{48}/V_8$ on
$w_2,\dots,w_8$ ($w_1$ is fixed) realize them as transitive
subgroups of $S_7$, so the results also follow from Theorem 1.1 (3).

\renewcommand{\refname}{\centering{References}}

\end{document}